\documentclass[12pt]{article}

\usepackage{graphicx}
\usepackage{color}
\usepackage{amsfonts}
\usepackage{amsmath}
\usepackage{natbib}
\usepackage{multirow}
\usepackage{amssymb}

\usepackage{graphicx}
\usepackage{amsfonts}
\usepackage{amsmath}
\usepackage[utf8]{inputenc}
\usepackage{url}
\usepackage{color}
\usepackage{mathtools}
\usepackage{multirow}
\usepackage{amsthm}

\begin{document}
\newcommand{\real}{\mathbb{R}}
\newcommand\numberthis{\stepcounter{equation}{}\tag{\theequation}}
\newcommand*{\medcup}{\mathbin{\scalebox{1.5}{\ensuremath{\cup}}}}%

\newcommand{\lip}{{\rm Lip}}

%\numberwithin{equation}{section}
%\theoremstyle{plain}
\newtheorem{theorem}{Theorem}[section]
\newtheorem{prop}{Proposition}[section]
\newtheorem{lema}{Lemma}[section]
\newtheorem{corolary}{Corollary}[section]
\newtheorem{remark}{Remark}[section]
\newtheorem{ex}{Example}[section]
\newtheorem{ass}{Assumption}

\title{\sc On approximate validation of models: A Kolmogorov-Smirnov based approach.\footnote{Research partially supported by 
FEDER, Spanish Ministerio de Econom\'{\i}a y Competitividad, grant MTM2017-86061-C2-1-P and Junta de Castilla y Le\'on, grants VA005P17 and VA002G18.}}

\author{Eustasio del Barrio,    Hristo Inouzhe and Carlos Matr\'an\\
\textit{Departamento de Estad\'{\i}stica e Investigaci\'on Operativa and IMUVA,}\\
\textit{Universidad de Valladolid. SPAIN} }
%\date{}
\maketitle

\begin{abstract}
 Classical tests of fit typically reject  a model for large enough real data samples. In contrast, often in statistical practice a model offers a good description of
the data even though it is not the ``true'' random generator. We consider a
more flexible approach based on contamination neighbourhoods around a model.
Using trimming methods and the Kolmogorov metric we introduce a functional statistic  measuring departures from a  contaminated model and the associated estimator corresponding to its sample version. We show how this estimator allows testing of fit
for  the (slightly) contaminated model vs sensible deviations from it, with uniformly exponentially small
type I and type II error probabilities. We also address the asymptotic behavior of the estimator showing that, under suitable regularity conditions,  it asymptotically behaves as   the supremum of a Gaussian process. As an application we explore methods of 
comparison between descriptive models based on the paradigm of model falseness. We also include some connections of our approach with the False-Discovery-Rate setting, showing competitive behavior when estimating the contamination level, although applicable in a wider framework.
\end{abstract}


\section{Introduction}

Classical Goodness of Fit  tests try to establish if there is enough statistical evidence to reject the null hypothesis,
which usually is a fixed generating mechanism. These procedures behave fairly well for moderate data sizes, but can become
excessively rigid in the presence of large sample sizes. This fact was already noted   for the chi-squared statistic 
in \cite{berkson1938some}  and interpreted by many authors as an indication of model falseness 
leading to statements such as `for every data generating mechanism there exists a sample size at which the model failure will become obvious' 
(see  \cite{modelAdec}). The issue has been approached in different ways  (see e.g. \cite{hodges1954testing},
\cite{simTrim}, \cite{Munk},\dots),  sharing the  idea that we should broaden the null hypothesis  to include
useful nearby models. Usually this is also accompanied by a gain in robustness in the new proposals.

However, considering the celebrated Box's phrase `essentially, all models are wrong, but some
are useful', even under the paradigm of model falseness, rejecting a model would  not be a satisfactory goal. If 
all models are false, and at a certain point, with enough data, we are able to reject the model, we could provide some measure 
of how useful  of how good  is compared to other models.  	This topic
is addressed in  \cite{davies1995data,Davies2016} from the perspective that a useful model is  anyone able to 
generate   similar samples  to the available data. 
Let us present our framework to revisit both topics from a novel point of view.

Often, some feature of a predominant population is clearly different from that of another minority population, simply because of its different eating or cultural habits. In either of these situations, a data sample of that feature taken from the general population will include data that do not come from and do not look like those arising from the predominant one. Consequently,  the statistical inference on the main population should be made taking into account the presence of atypical data. As a first ingredient, to address this goal, we resort to a suggestive model introduced  in  \cite{huber1964robust}, becoming one of the very basis of  Robust Statistics: An ($\alpha$-)contamination neighbourhood (CN) of a probability distribution $P_0$ is the set 
of probability distributions
\begin{equation}
\label{def_cont_neigh}
\mathcal{V}_\alpha(P_0)=\{(1-\alpha)P_0+\alpha Q: Q\in\mathcal{P}\},
\end{equation}
where $\mathcal{P}$ is the set of all probability distributions in the space (throughout the paper  the real 
line $\mathbb{R}$).
For a given probability $P_0$ and a particular value $\alpha_0\in [0,1)$, a probability $P$ in $\mathcal{V}_{\alpha_0}(P_0)$ would generate samples with an approximate   $(1-\alpha_0)\times100\%$   of data coming from $P_0$. Also we must note the  use of particular contamination models in different statistical problems, stressing its role on the False-Discovery-Rate  (FDR) setting (as considered e.g. in \cite{Genovese}). We briefly comment on the relation of our approach with that in Section \ref{FDR}. Of course, if an `outlying label' were available for the data coming from  the contaminating  distribution, $Q$, removing the labeled data would produce a legitimate sample from $P_0$. The relevant fact is that CN's are  related to trimmings (see \cite{optTrim}) by
\begin{equation}\label{equivalencia}
P\in \mathcal{V}_\alpha(P_0) \Longleftrightarrow P_0\in R_\alpha(P),
\end{equation}
 where $R_\alpha(P)$ denotes the  set of $\alpha$-trimmings of the probability distribution $P$,
\begin{equation}
\label{def_alpha_trimm}
R_\alpha(P)=\big\{Q\in\mathcal{P}:Q\ll P,\,{\textstyle \frac{dQ}{dP}\leq\frac{1}{1-\alpha}}\, P\text{-a.s.}\big\}.
\end{equation}
This means that {\bf an $\alpha$-trimming, $Q$, of $P$} is characterized by a down-weighting function $f$ such that 
$0\leq f \leq 1$ and $Q(B)=\frac 1 {1-\alpha}\int_Bf(x)P(dx)$ for all measurable sets $B$ in $\mathbb R$. In contrast with 
the hard 0-1 (trimmed/non-trimmed)  trimming practice in data analysis, this concept allows for gradually diminishing/enhancing
the importance of points  in the sample space. 
Relation (\ref{equivalencia}) allows us to work with trimmings, instead of CN's, taking 
advantage of the underlying meaning of trimming and its mathematical properties.     If $F$ and $F_0$ are distribution functions 
(d.f.'s in the sequel), we will also use $R_\alpha(F)$ and $\mathcal{V}_\alpha(F_0)$, with the same meanings as before, but 
defined in terms of d.f.'s.  

The natural absence of an outlying label  has been traditionally substituted by more or less orthodox trimming criteria, including the oldest consisting in trimming just the extreme values, carrying out the analysis with the remaining data. Recently, mainly in conection with two-sample problems (see e.g. \cite{TrimmedCompar, optTrim,simTrim,Alvarez-Esteban2016}), optimal trimmings have been introduced as the nearest ones to the original model, according to some probability distance or dissimilarity measure. This role will be played here by the Kolmogorov (or $L_\infty$-)distance between 
d.f.'s on the line, namely, 
$$d_K(F,G)=\sup_{x\in \mathbb R}|F(x)-G(x)|,$$ 
(we will often use the notation $\|F-G\|$ for $d_K(F,G)$).

{In this work,  we  develop a robust hypothesis 
testing procedure based on the previous considera\-tions. Moreover, under the paradigm of a false-model world, we use the elements involved in the procedure
to suggest some tools for comparing models or to determining the usefulness of   particular models.}

The use of CN's, through their connection with trimmings, leads to consider 
$\mathcal{V}_\alpha(F_0)$ to be the `reasonable' model. Notice that (see Example \ref{EjemploUniformes}), this approach differs from that based just on $d_K$-neighbourhoods  of $F_0$,  which would have a different meaning       (see \cite{Owen} for this and other classic approaches).
As  relation (\ref{equivalencia2}) shows, (\ref{equivalencia}) is also equivalent to $d_K(F_0,R_\alpha(F))=0$,  giving to the  `trimmed Kolmogorov distance' functional 
\begin{equation}\label{statfunc}
d_K(F_0,R_\alpha(F)):=\min_{\tilde{F}\in R_\alpha(F)}d_K(F_0,\tilde{F}),
\end{equation}
and to the  plug-in estimator  $d_K(F_0,R_\alpha(F_n))$,
a main role into our analysis. (Here $F_n$ is the empirical d.f. based on  a sample  of $n$ independent
random variables with common d.f. $F$).  In particular, we address the possibilities of testing $H_0:d_K(F_0,R_\alpha(F))=0 \text{ vs. } H_1:d_K(F_0,R_\alpha(F))>0,$ where `reasonable'  is controlled by the trimming level $\alpha$.
Related null hypotheses have already been considered making use of 
different probability metrics or different neighbourhoods.  In \cite{optTrim,simTrim}, the $L_2$-Wasserstein distance is used in a two-sample version. Previous approaches based on particular 
trimming procedures were considered in \cite{Munk} and \cite{TrimmedCompar}. The Kolmogorov-Smirnov test is probably the most widely 
used goodness of fit test, therefore the $d_K$-metric provides a privileged setting to develop our approach.  Notice that, in \cite{Hristo}, we have included most of the mathematical tools involved in this problem. This includes     existence and characterization of (a particular)  minimizer, and even a result on directional differentiability, which will be used here. 
%Under the present framework we will be able to obtain exact exponential bounds for error probabilities. This last feature is highly desirable for any hypothesis 
%testing procedure. 

As shown in \cite{uecBa}, for any 
distance $d$  dominating the total variation distance, testing the null hypothesis $P=P_0$ vs. the alternative $d(P,P_0)\geq \rho\  (>0)$, 
makes  generally unachievable to get exponential bounds for the involved errors. The test provided in Section 3 has exponentially small error probabilities for testing the null $H_0:d_K(F_0,R_\alpha(F))=0$ (equivalently, 
$H_0: F\in \mathcal{V}_\alpha(P_0)$) against the alternative $d_K(F_0,R_\alpha(F))>\rho$. The test is uniformly consistent
(type I and type II error probabilities tend to 0 uniformly) for detecting alternatives $d_K(F_0,R_\alpha(F))>\eta_n/\sqrt{n}$ with $\eta_n/\sqrt{n}\to 0$  
if $\eta_n\to\infty$.

Also, in Section  \ref{CLT}, we provide asymptotic  theory for $d_K(F_0,R_\alpha(F_n))$ for inferential purposes. It includes an extension of Theorem 2 in \cite{dkR}  for flexible null hypotheses.

The second main goal in this paper is to provide tools to compare  different models when the null hypothesis is rejected.
Under the model falseness paradigm, 
\cite{davies1995data,Davies2016} introduce the idea of adequacy region (for a data set) as the set of probabilities 
in a model whose samples would typically look like the actual data. Also
\cite{rudas1994new} proposes the very natural concept  of index of fit, namely, the contamination level 
necessary to make the random generator of the data a contaminated member of the model. 
The proposal in \cite{rudas1994new}, as well as its modification in \cite{parMult},  deal  with multinomial models. In our setup 
we  consider the trimmed Kolmogorov (tK) index of fit, $\alpha^*$, defined by
\begin{equation}
\label{eq_emp_alpha}
\alpha^*=\min\{\alpha:\ d_K(F_0, R_\alpha(F))=0\}.
\end{equation}
This is the minimum contamination level $\alpha$ for which $F$ is a contaminated version  of $F_0$. 
This works in a very general 
setup, since we impose no constraints on $F$ and $F_0$.  This is in contrast with  the methodology involved in the control of FDR, which takes advantage of the  dominated contamination model. With the methodology developed here,  it is fairly easy to 
calculate the empirical version of $\alpha^*$ for a particular data set.
Using our asymptotic theory for $d_K(F_0,R_\alpha(F_n))$ we propose a consistent estimator for $\alpha^*$ in Section 4. We also provide  comparisons with some methodologies developed in the FDR setting (as considered in \cite{Meinshausen}) for estimating the proportion of false null hypotheses. 

A related approach for comparing the quality of different models to describe the data  is based on credibility indices, as introduced in \cite{modelAdec}. Given a  goodness of fit procedure, the credibility index allows comparison between models based on the minimal sample size $n^*$ for which subsamples of size $n^*$ of the original data (of size $n$) 
reject the null hypothesis 50\% of times. The idea behind this index is that for large samples,  goodness of fit tests will  
very likely  reject the null hypothesis, while often for smaller sub-samples  the null would not be rejected. 
Of course, these credibility indices have to be estimated from the data. The proposal in \cite{modelAdec} is to use subsampling to perform
this estimation. However, the accuracy of the subsampling approximation is limited to small (as compared to the complete sample) subsample sizes.
Here we show how our asymptotic theory for $d_K(F_0,R_\alpha(F_n))$ can provide further information about the credibility indices.

Summarizing, the paper addresses the analysis and applications of $d_K(F_0,R_\alpha(F))$, the `trimmed Kolmogorov distance'. Section 2 is devoted to collect the mathematical bases  and provide a fast algorithm for computation on sample data. The analysis of the proposed testing procedure is carried in  Section \ref{HIPTEST}. In Section \ref{Comparison} we show how to apply this test to
credibility analysis and develop some results about the tK-index of fit and the related acceptance regions. The basis for that approach relies on the CLT for the trimmed Kolmogorov distance (see Theorem \ref{tma_tcl}). Section \ref{FDR} includes some relations with the FDR setting and comparisons between several estimators of the contamination index $\alpha.$ In Section \ref{Simulations} 
we illustrate the previous techniques to compare descriptive models over simulated and real data examples. In the last section we 
briefly discuss the results. Finally, the proof of the main result in the paper, the CLT for the trimmed Kolmogorov distance, is given in  the Appendix.

\section{Trimming and Kolmogorov distance}

We keep the notation used in the Introduction and notice that the set $R_\alpha(F)$ can be also characterized, as showed in 
\cite{TrimmedCompar} (Proposition 2.2 in \cite{optTrim} gives a more general result), in terms of the set  of 
$\alpha$-trimmed versions of the uniform law $U(0,1)$. Let  $\mathcal{C}_\alpha$ be the set   of absolutely continuous functions $h:[0,1]\mapsto[0,1]$, such that
$h(0)=0, h(1)=1$, with derivative $h'$ verifying $0\leq h'\leq\frac{1}{1-\alpha}\,\text{a.e.}.$ Then, the composition of the functions  $h$ and $F$: $F_h=h\circ F$ gives the useful parameterization
\begin{equation}\label{calpha}
R_\alpha(F)=\{F_h: h \in \mathcal{C}_\alpha\}.
\end{equation}

The set $R_\alpha(F)$ is convex and also well behaved w.r.t.  weak convergence of probabilities and 
widely employed probability metrics (see Section 2 in \cite{optTrim}). 
As showed in \cite{Hristo}, $R_\alpha(F)$ keeps several nice properties under $d_K$; we include below the most relevant ones.
\begin{prop}
\label{dkcompact}
For $\alpha \in [0,1)$, if  $F$, $G$ with or without suffixes are d.f.'s:
\begin{itemize}
\item[(a)] $R_\alpha(F)$ is compact w.r.t. $d_K$.
\item[(b)] $d_K(F_0,R_\alpha(F))=\min_{\tilde{F}\in R_\alpha(F)}\|\tilde{F} -F_0\|=\min_{h\in \mathcal{C}_\alpha}\|h\circ F -F_0\|$.
\item[(c)] $
|d_K(G_1,R_\alpha(F_1))-d_K(G_2,R_\alpha(F_2))| \leq d_K(G_1,G_2)+ {\textstyle \frac 1 {1-\alpha}}{d_K(F_1,F_2)}.
$
\item[(e)] If  $d_K(F_n,F)\to 0,$ then:

\item[e1)] for every $\tilde{F}\in R_\alpha(F),$  there exist $\tilde{F}_n\in R_\alpha(F_n), n\in \mathbb N$ such that $d_K(\tilde{F}_n,\tilde{F})\to 0.$
\item[e2)] if $\tilde{F}_n \in R_\alpha(F_n), n\geq 1$, then there exists some $d_K$-convergent subsequence $\{\tilde{F}_{n_k}\}$. If $\tilde{F}$ is the limit of 
such a subsequence, necessarily $\tilde{F}\in  R_\alpha(F)$.
\item[e3)] if, additionally,  $\{G_n\}$ is any sequence of d.f.'s  such that $d_K(G_n,G)\to 0,$ then $d_K(G_n,$ $R_\alpha(F_m)) \to d_K(G,R_\alpha(F))$ as $n,m \to \infty.$
\end{itemize}
\end{prop}

Immediate consequences of Proposition \ref{dkcompact} are that for $\alpha \in [0,1)$: 
\begin{equation}\label{minimizers}
\mbox{There exists } \ \tilde F_0\in  R_\alpha(F) \ \mbox{ such that } \ d_K(F_0,\tilde F_0)= d_K(F_0,R_\alpha(F)), \ \mbox{ and } 
\end{equation}
\begin{equation}\label{equivalencia2}
F \in \mathcal{V}_\alpha(F_0) \ \mbox{ if and only if } d_K(F_0,R_\alpha(F))=0.
\end{equation}
Moreover, by convexity of $R_\alpha(F)$, the set of optimally trimmed versions of $F$ associated to problem (\ref{minimizers}) is also convex.
However, guarantying uniqueness of the minimizer (as it holds w.r.t. $L_2$- Wasserstein metric by Corollary 2.10 in \cite{optTrim}) is not possible. 
Mention apart, by its statistical interest, merits the  the following  consistency result, which is straightforward from Glivenko-Cantelli theorem and item e3) above.

\begin{prop}[Consistency of trimmed Kolmogorov distance]
\label{consistency}
Let $\alpha \in [0,1)$ and  $\{F_n\}$ be the sequence of empirical d.f.'s based on a 
sequence $\{X_n\}$ of independent random variables with distribution function $F$. If $\{G_n\}$ is any  sequence of distribution 
functions $d_K$-approximating the d.f. $G$ (i.e. $d_K(G_n,G)\to 0$), then:
$$d_K(G_n,R_\alpha(F_m))\to d_K(G,R_\alpha(F)),\ \mbox{ as } n,m \to \infty, \ \mbox{ with probability one.}$$
\end{prop}

While in other contexts the roles played by discarding contamination (by  trimming)  and the distance under consideration seem to be clear, here the nature of Kolmogorov distance can lead to a distorted picture. To give some light on these roles, we include  a very simple example based on uniform laws that allows explicit  computations. We also must note that (as commented in \cite{simTrim}) contamination neighbourhoods have been extended in several ways; notably Rieder's neighborhoods of a probability comprise contamination as well as total variation norm neighborhoods.

\begin{ex}\label{EjemploUniformes} {\em Contamination vs  $d_K$-based neighbourhoods.} { Let us fix $F_0$ to be the $U(0,1)$ d.f. and consider the following scenarios for $F$
\begin{itemize}
\item[ i)] $F$ the d.f. of an $U(0,1+\varepsilon)$ or an $U(-\varepsilon,1)$ law. Then $d_K(F_0,F)=\frac{\varepsilon}{1+\varepsilon}$ and $d_K(F_0,R_\alpha(F))=\frac{\varepsilon-\alpha}{1+(\varepsilon-\alpha)}$ if $0\leq\alpha\leq \varepsilon$ (and 0 if $\alpha \geq \varepsilon).$
\item[ii)] $F$ the d.f. of a $U(0,1-\varepsilon)$ law. Then $d_K(F_0,F)=\varepsilon$ and $d_K(F_0,R_\alpha(F))=\varepsilon$ for every $0\leq \alpha< 1$.
\end{itemize}
In fact, the first situation involves a contamination of exact size $\varepsilon$ of $F_0$, because $F=(1-\varepsilon)F_0+\varepsilon F'$ where $F'$ is the d.f. of an $U(1,1+\varepsilon)$ or an $U(-\varepsilon,0)$ law. In contrast, the second one does not fit in the contamination model at all. The following scenario includes inner contamination at the support of $F_0$, adding some complexity to the analysis:
\begin{itemize}
\item[ iii)] $F=(1-\varepsilon)F_0+\varepsilon F'$, where $F'$ is the d.f. of a $U(a,b)$ law with $0<a<b<1$. Then $d_K(F_0,F)=\varepsilon \sup\{a,1-b\},$ and for $0\leq \alpha \leq \varepsilon$: $d_K(F_0,R_\alpha(F))=(\varepsilon-\alpha)\sup\{a,1-b\},$ if $0<a<b \leq1/2$ else $1/2\leq a<b<1$. If $0<a\leq 1/2 <b<1$, then for $0<\alpha<\varepsilon_0:=\varepsilon \frac{|a+b-1|}{b-a}$, we would have $d_K(F_0,R_\alpha(F))=(\varepsilon-\alpha)\sup\{a,1-b\},$ while for $\varepsilon_0 \leq \alpha \leq \varepsilon,$  defining $\gamma=|1/2-\sup\{a,1-b\}|$, we would have $d_K(F_0,R_\alpha(F))=[1/2-\gamma(\varepsilon-\alpha)/(\varepsilon-\varepsilon_0)](\varepsilon-\alpha)$.
\end{itemize}
The analysis above shows that the effect of optimal trimming according to the $d_K$-distance strongly depends on several factors. Notably, they include the presence or not of a contaminating part, but also its  spread and relative position. 
\hfill $\Box$
}
\end{ex}

Throughout this paper we make frequent use of the quantile function. Given a d. f. $F$, we write $F^{-1}$ for the associated quantile function. Recall that it is just the left-continuous inverse  of the d.f. $F$, namely, 
$F^{-1}(t):=\inf\{x\ | \ t\leq F(x)\}$. It allows a useful representation  of  the corresponding distribution because, if $U$ is a uniformly distributed 
$U(0,1)$ random variable,  $F^{-1}(U)$ has d.f. $F$. Moreover, if $X$ has a continuous d.f. $F$, $F_0\circ F^{-1}$ is  easily seen to be
the quantile function associated to
the r.v. $Y=F_0(X)$. As we showed in \cite{Hristo}, under some regularity assumptions, $d_K(F_0,R_\alpha(F))$ can be expressed in terms of the function
$F_0\circ F^{-1}$.    This fact allows the practical computation of  $d_K(F_0,R_\alpha(F_n))$ 
when $F_n$ is an empirical d.f. based on a data sample $x_1,\dots,x_n$, and even that of $d_K(F_0,R_\alpha(F))$ for theoretical distributions (see Example \ref{EjemploNormal1}). For the sake of completeness, we include below these results and a theorem which is a fundamental tool for our goals. It gives an explicit characterization of a solution of the corresponding optimization problem (see Theorem 2.5 in \cite{Hristo} for a proof). 

\begin{lema} \label{ref1}
If $F, F_0$ are continuous d.f.'s  and $F$ is additionally strictly increasing then
$$d_K(F_0,R_\alpha(F))=\min_{h\in \mathcal{C}_\alpha}\|h-F_0\circ F^{-1}\| \
\mbox{ and } \
d_K(F_0, R_\alpha(F_n))=\min_{h\in \mathcal{C}_\alpha}\|h-F_0 \circ F_n^{-1}\|.$$
\end{lema}

\begin{theorem}
\label{prop_disting_h}
Assume $\Gamma:[0,1]\to [0,1]$ is a continuous nondecreasing function. Define $G(t)=\Gamma (t)-\frac t {1-\alpha}$, $U(t)= \sup_{t\leq s\leq 1} G(s)$, $L(t)=\inf_{0\leq s\leq t}G(s)$ and
$$\tilde{h}_\alpha(t)=\max\left(\min\left( {\textstyle \frac{U(t)+L(t)}{2} }, 0\right),{\textstyle  \frac{-\alpha}{1-\alpha}}\right).$$
Then $h_\alpha:=\tilde{h}_\alpha+\frac{\cdot}{1-\alpha}$ is an element of $\mathcal{C}_\alpha$, and
$$\min_{h\in \mathcal{C}_\alpha}\|h-\Gamma\|=\|h_\alpha-\Gamma\|=\|\tilde{h}_\alpha-G\|.$$
\end{theorem}

 Note that the assumption on $\Gamma$  is always verified when $\Gamma=F_0\circ F^{-1}$, and that taking  right and left limits
at 0 and 1, respectively, we can assume that $F_0\circ F^{-1}$ is a nondecreasing (and left continuous) function from $[0,1]$ to $[0,1]$.

 A key aspect in Theorem \ref{prop_disting_h} is that, although not necessarily unique,  $h_\alpha$ is an optimal trimming function in the sense described above. However, from the point of view of asymptotic theory, Theorem \ref{prop_disting_h} is the key to our 
Theorem \ref{tma_tcl} in Section \ref{Comparison}. Moreover, from a practical point of view, it yields a simple algorithm for the computation of $d_K(F_0,R_\alpha(F_n))$, as follows. 

Assume $X_1,\ldots,X_n$ are 
i.i.d. observations from the continuous and strictly increasing d.f. $F$ and assume that $F_0$ is continuous. From Lemma \ref{ref1} and Theorem \ref{prop_disting_h} we know that 
$d_K(F_0,R_\alpha(F_n))=\|\tilde{h}_{\alpha,n}-G_n\|$, where $G_n(t)=H_n^{-1}(t)-\frac{t}{1-\alpha}$, $H_n^{-1}$ is the empirical quantile function of the transformed data, $Y_i=F_0(X_i)$,
$U_n(t)= \sup_{t\leq s\leq 1}G_n(s)$, $L_n(t)=\min_{0\leq s\leq t}G_n(s)$ and $$\tilde{h}_{\alpha,n}(t)=\max\left(\min\left( \frac{U_n(t)+L_n(t)}{2} , 0\right), \frac{-\alpha}{1-\alpha}\right).$$
Denote by $Y_{(1)}\leq\cdots\leq Y_{(n)}$ the ordered (transformed) sample. Note that $G_n(t)=Y_{(i)}-\frac t{1-\alpha}$ if $t\in (\frac{i-1}n,\frac i n]$, while
$\tilde{h}_{\alpha,n}$ is a nonincreasing function and this implies that
$$\|\tilde{h}_{\alpha,n}-G_n\|=\max_{1\leq i\leq n}\left(\max(G_n({\textstyle \frac {i-1} n}+)-\tilde{h}_{\alpha,n}({\textstyle \frac {i-1} n}),
\tilde{h}_{\alpha,n}({\textstyle \frac {i} n})-G_n({\textstyle \frac {i} n}) ) \right),$$
with $G_n({\textstyle \frac {i-1} n}+)=Y_{(i)}-\frac{i-1}{n(1-\alpha)}$, $G_n({\textstyle \frac {i} n})=Y_{(i)}-\frac{i}{n(1-\alpha)}$. 
For the computation of $\tilde{h}_{\alpha,n}({\textstyle \frac {i} n})$ we note that $U_n(\frac i n)=\max_{i\leq j\leq n-1} G_n({\textstyle \frac {j} n}+)$ and 
$L_n(\frac i n)=\min_{1\leq j\leq i} G_n({\textstyle \frac {j} n})$ for $i=1,\ldots,n-1$. Summarizing, we see that $d_K(F_0,R_\alpha(F_n))$ can be computed through the following algorithm.

\bigskip
\textit{\bf Algorithm for the computation of $d_K(F_0,R_\alpha(F_n))$:}

\begin{itemize}
 \item compute $Y_i=F_0(X_i)$, $i=1,\ldots,n$; sort $Y_{(1)}\leq \cdots \leq Y_{(n)}$.
 \item compute $g_{i+}=Y_{(i+1)}-\frac{i}{n(1-\alpha)}$, $i=0,\ldots, n-1$; $g_{i-}=Y_{(i)}-\frac{i}{n(1-\alpha)}$, $i=1,\ldots,n$.
 \item compute $u_i=\max_{i\leq j\leq n-1} g_{j+}$, $l_{i}=\min_{1\leq j\leq i} g_{j-}$, $i=1,\ldots,n-1$.
 \item set $h_0=0$, $h_n=-\frac \alpha{1-\alpha}$ and $h_i=\max(\min(0,\frac {u_i+l_i}2),-\frac \alpha{1-\alpha})$, $i=1,\ldots,n-1$.
 \item compute 
 $$d_K(F_0,R_\alpha(F_n))=\max_{1\leq i\leq n}\big(\max(g_{(i-1)+}-h_{i-1}, h_i-g_{i,-}\big).$$
\end{itemize}

Beyond this algorithm for the empirical case, Theorem \ref{prop_disting_h} provides a simple way for the computation
of theoretical trimmed Kolmorogov distances. Example 2.1 in \cite{Hristo} analyzes the problem in Gaussian model. Let us include here a summary for  illustration of this use.

\begin{ex}\label{EjemploNormal1} {\em Trimmed Kolmogorov distances in the Gaussian model.} {
Consider the case $F_0=\Phi$, $F=\Phi((\cdot-\mu)/\sigma)$, where $\Phi$ denotes the standard normal d.f., $\mu\in\mathbb{R}$ and $\sigma>0$.
Here we have $H^{-1}(t):=F_0\circ F^{-1}(t)=\Phi(\mu + \sigma\Phi^{-1}(t))$. 
We focus on the cases $\sigma=1$, $\mu\ne 0$ and $\mu=0$, $\sigma\ne 1$ (see \cite{Hristo} for details).

\smallskip
If $\sigma=1$ and $\mu\neq0$ then
\begin{equation}\label{normaldominated}
d_K(R_\alpha (N(\mu, 1)), N(0,1))= 
{\textstyle \Phi\big(\frac {|\mu|} 2  +\frac 1 {|\mu|} \log(1-\alpha)\big)-
\frac 1 {1-\alpha}\Phi\big(-\frac {|\mu|} 2+\frac 1 {|\mu|} \log(1-\alpha)\big)}.
\end{equation}

\smallskip
In the case $\mu=0$: 
$${ d_K(R_\alpha(N(0, \sigma^2)), N(0,1))=\begin{cases}{\textstyle\Phi\Big(\frac {-\sigma \frac \Delta 2 }{1-\sigma^2} \Big)-\frac 1 {1-\alpha}
\Phi\Big(\frac {- \frac \Delta 2 }{1-\sigma^2} \Big)},\quad \mbox{if }\sigma<1\\

0, \quad \quad \quad \quad \quad \quad \quad \quad \mbox{ if } 1\leq \sigma \leq 1/(1-\alpha)\\

{\textstyle\Phi\Big(\frac {\sigma \frac \Delta 2 }{\sigma^2-1} \Big)-\frac {\Phi\Big(\frac { \frac \Delta 2 }{\sigma^2-1} \Big)-\frac \alpha 2}{1-\alpha}, \mbox{ if } \sigma > 1/(1-\alpha)}\\

\end{cases}}$$
\hfill $\Box$
}
\end{ex}

\begin{figure}[htb]
\begin{center}
\includegraphics[scale=0.4]{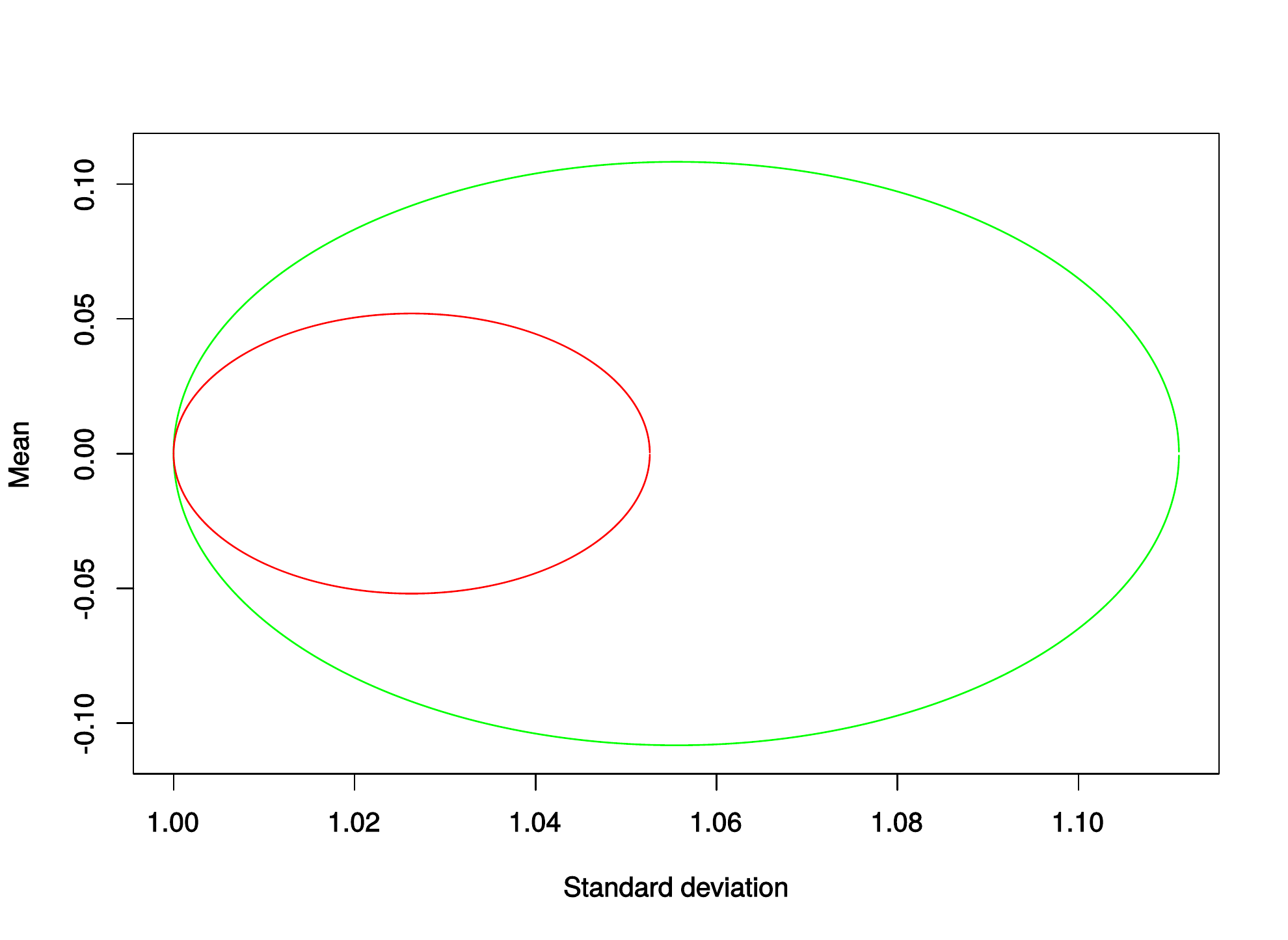}
\caption{Plot of regions containing the parameters compatible with $\alpha$-contamination 
neighbourhoods of $F_0\sim N(0,1)$,  for $\alpha=0.05$ (red) and $\alpha=0.1$ (green)}
\label{fig_accept_region}
\end{center}
\end{figure}

Relations (\ref{minimizers}) and (\ref{equivalencia2}) state the  link between CN's and  
trimming, opening ways to approximately validating a model making use of trimming through the Kolmogorov distance. We end this 
section showing how CN's and approximate validation in a parametric model setting 
can be related. For that task we focus on what are the parameters in the model leading to distributions in 
$\mathcal{V}_{\alpha}(F_0)$. As pointed out in \cite{davies1995data}, we should just consider models able to generate 
data similar to our sample. Moreover, distributions in a CN have an intuitive appeal and, 
if $\alpha$ is small, we can expect to be handling reasonable models. For instance, if $F_0 \sim N(0, 1)$ then we can 
calculate the  tolerance region given by the subset of normal distributions belonging to $\mathcal{V}_{\alpha}(F_0)$ in an elementary 
fashion. This provides an approximate picture of the kind of distributions present in the CN of $F_0$. 
These tolerance regions for $\alpha = 0.05$ and $\alpha = 0.1$ are shown in Figure \ref{fig_accept_region}.
Every combination of $(\tilde{\mu}, \tilde{\sigma})$ inside the green border is a normal distribution that belongs to
$\mathcal{V}_{0.1}(N(0,1))$. The same is true for the red border and $\mathcal{V}_{0.05}(N(0,1))$.

\section{Hypothesis testing}\label{HIPTEST}
 
To develop our approach for a testing procedure, throughout, $X_1,\dots,X_n$ will be  independent random variables with common  d.f. $F$, and  $F_n$ will be the corresponding empirical d.f. 
The main result, following the principles in \cite{uecBa}, concerns control of error probabilities: a test is uniformly consistent (UC) if 
both type I and type II error  probabilities (EI and EII in the sequel) converge uniformly to 0 as the sample size, $n\to \infty$, and 
it is uniformly exponentially consistent (UEC) if the error probabilities  are uniformly bounded by $e^{-rn}$ for  large $n$ and some $r>0$. To stress on the necessity of considering some separating zone between the null and the alternative, we include this previous slightly more general result.
\begin{prop}
\label{tma_test}
Given $0\leq \rho_1<\rho_2$, for testing
$H_0:d_K(F_0,R_\alpha(F))\leq\rho_1$  vs.  $H_1:d_K(F_0,R_\alpha(F))>\rho_2$, for every $0<\lambda< 1$ rejecting the null hypothesis when
$d_K(F_0, R_\alpha(F_n))>(1-\lambda)\rho_1+\lambda\rho_2$
is an uniformly exponentially consistent (UEC) test.
\end{prop}
\medskip
\noindent \textbf{Proof.}
  From Proposition 2.1 c) in \cite{Hristo}, we have the inequality $|d_K(F_0,R_\alpha(F_1))-d_K(F_0,R_\alpha(F_2))| \leq {\textstyle \frac 1 {1-\alpha}}\|F_1-F_2\|$, thus  for EI: 
\begin{align*}
&P_{F}\Big(d_K(F_0,R_\alpha(F_n))>(1-\lambda)\rho_1+\lambda\rho_2\Big)\\
&\leq P\Big(d_K(F_0,R_\alpha(F_n))-d_K(F_0,R_\alpha(F))>\lambda\rho_2-\lambda\rho_1\Big)\\
%&\leq P\Big(d_K(R_\alpha(F_n),R_\alpha(F))>\frac{\rho_2-\rho_1}{2}\Big)\\
&\leq P\Big(\frac{1}{1-\alpha}\sup_x|F_n(x)-F(x)|>\lambda(\rho_2-\rho_1)\Big)\\
&=P\Big(\sup_x\sqrt{n}|F_n(x)-F(x)|>\sqrt{n}(1-\alpha)\lambda(\rho_2-\rho_1)\Big)\\
&\leq 2e^{-2\lambda^2 n(1-\alpha)^2(\rho_2-\rho_1)^2}.\numberthis\label{EI}
\end{align*}
Note that the last bound follows from the \cite{Massart} version of the Dvoretsky-Kiefer-Wolfowitz inequality.

To handle EII (thus if $d_K(F_0,R_\alpha(F))> \rho_2$), we have
\begin{flalign*}
&P_{F}\Big(d_K(F_0,R_\alpha(F_n))\leq(1-\lambda)\rho_1+\lambda\rho_2\Big)&\\
&\hspace*{50pt}=P_{F}\Big(\rho_2-d_K(F_0,R_\alpha(F_n))\geq\rho_2-((1-\lambda)\rho_1+\lambda\rho_2)\Big)&\\
&\hspace*{50pt}\leq P\Big(d_K(F_0,R_\alpha(F))-d_K(F_0,R_\alpha(F_n))>(1-\lambda)(\rho_2-\rho_1)\Big)&\\
%&\leq P\Big(d_K(R_\alpha(F),R_\alpha(F_n))>\frac{\rho_2-\rho_1}{2}\Big)\\
&\hspace*{50pt}\leq P\Big(\sup_x\sqrt{n}|F_n(x)-F(x)|>\sqrt{n}(1-\alpha)(1-\lambda)(\rho_2-\rho_1)\Big)&\\
&\hspace*{50pt}\leq 2e^{-2(1-\lambda)^2 n(1-\alpha)^2(\rho_2-\rho_1)^2}.&\numberthis\label{EII}
\end{flalign*}
\quad $\Box$
\vspace{5mm}

As an easy consequence, taking $\rho_1=0$ and $\rho=\rho_2$, we get:
\begin{theorem}\label{Testing}
Given $\rho>0$, for testing
\begin{equation}
\label{test0}
H_0:d_K(F_0,R_\alpha(F))=0\quad\text{ vs. }\quad H_1:d_K(F_0,R_\alpha(F))>\rho,
\end{equation} for every $0<\lambda< 1$ the critical region
$d_K(F_0, R_\alpha(F_n))>\lambda\rho$ defines
 an uniformly exponentially consistent (UEC) test.
\end{theorem}

Since the null  hypothesis includes all the contamination versions (of $\alpha$-level) of $F_0$,   rejection means that the generator of the sample is far enough of   any such a contaminated version. Theorem \ref{Testing} guarantees that alternatives will be quickly detected when farness is measured through the $d_K$-distance.  

In statistical practice, it could be wiser to change the alternative hypothesis and 
make it sample size dependent.
That leads to consider tests of the form
\begin{equation}
\label{test1}
H_{0,n}:d_K(F_0,R_\alpha(F))=0\quad\text{ vs. }\quad H_{1,n}:d_K(F_0,R_\alpha(F))>\rho_n,
\end{equation}
for $\rho_n=\rho(n)>0$, and rejection when $d_K(F_0,R_\alpha(F_n))>\lambda\rho_n$. For instance, taking
$\rho_n = \eta_n/\sqrt{n}\to 0$, and $\eta_n\to \infty$ results in an uniformly consistent test. Uniform consistency is weaker  than 
uniform exponential consistency, but it allows to detect, for example, alternatives at a distance $\log(n)/\sqrt{n}$. Also,
we can consider $\lambda$ as a tuning parameter which can help if we have some additional information or if we want more or less conservative tests with respect to 
EI and EII probabilities (of course, when $\rho=0$, $\rho_1=\rho_2$ or $\lambda=0$ or $\lambda=1$, some bounds 
are meaningless and we can not assure uniform consistency with the previous procedure). Alternatively, 
we may look for the smallest possible values for $\rho_n$, while still controlling  EI and EII.
From (\ref{EI}) and (\ref{EII}) note that if $\rho_n$ is $o(n^{-1/2})$ we would lose the control of the errors, since $n\rho_n^2\to 0$ as $n\to\infty$. This leads us to choose $\rho_n$ as $O(n^{-1/2})$, or, fixing some  $\rho>0$:
\begin{equation}
\label{rhon}
\rho_n=\frac{\rho}{\sqrt{n}}
\end{equation}
 {Now,  if we fix $0<\epsilon_1,\,\epsilon_2<1$, looking for a rejection threshold, $\lambda \rho_n$,  for which
$$EI\leq \epsilon_1\quad\text{and}\quad EII\leq\epsilon_2,$$
 we get
$2e^{-2\lambda^2 (1-\alpha)^2\rho^2}=\epsilon_1$ and $\epsilon_1e^{4\lambda-2}=\epsilon_2.$
With a bit of algebra we get
\begin{equation}\label{rhon_explicitb}
\rho=\frac{1}{(1-\alpha)\lambda}\sqrt{\frac{1}{2}\log\frac{2}{\epsilon_1}},\quad
\lambda=\frac{1}{2}+\frac{1}{4}\log\frac{\epsilon_2}{\epsilon_1},
\end{equation}
imposing  $\epsilon_1e^{-2}<\epsilon_2<\epsilon_1 e^2$, which gives the optimal boundary level
\begin{equation}
\label{rhon_explicit}
\rho_n=\frac{\rho}{\sqrt{n}}=\frac{1}{(1-\alpha)\lambda}\sqrt{\frac{1}{2n}\log\frac{2}{\epsilon_1}}.
\end{equation}
Relations (\ref{rhon_explicitb}) and (\ref{rhon_explicit}) summarize the balance among the different elements. 
Ideally, we look for small $\rho_n$, $\epsilon_1$ and $\epsilon_2$ but, paying the price for our demands, $\rho_n$ grows as $\epsilon_1$ gets 
smaller and as $\epsilon_2$ gets more similar to $\epsilon_1$. Therefore, we need to make sensible choices for $\epsilon_1$ and $\epsilon_2$. 
In Table \ref{table_rhon_behaviour} we show some examples of the  mentioned behaviour. For instance, fixing $\epsilon_1=0.01$ 
and $\epsilon_2=0.05$ seems a sensible choice, giving a fairly low $\rho_{1000}$ while keeping low error probabilities. 

\begin{table}[ht]
\centering
\caption{Values associated to error bounds for $\alpha = 0.1$ and $N=1000$.}
\label{table_rhon_behaviour}
%\begin{adjustbox}
\small
\begin{tabular}{|c|c|c|c|c|c|c|c|c|c|c|c|}
\hline
{ EI}  & {EII}  & {$\lambda$} & {$\rho_{1000}$} & { EI}  & {EII}  & {$\lambda$} & {$\rho_{1000}$} & { EI}  & {EII}  & {$\lambda$} & {$\rho_{1000}$} \\ \hline
0.1 & 0.5  & 0.90    & 0.048      & 0.05 & 0.25 & 0.90    & 0.053      & 0.01 & 0.05  & 0.90    & 0.063   \\ \hline
0.1 & 0.1  & 0.50    & 0.086      & 0.05 & 0.05 & 0.50    & 0.095      & 0.01 & 0.01  & 0.50    & 0.114   \\ \hline
0.1 & 0.02 & 0.10    & 0.440      & 0.05 & 0.01 & 0.10    & 0.489      & 0.01 & 0.002 & 0.10    & 0.586   \\ \hline
\end{tabular}
%\end{adjustbox}
\end{table}

An appealing goal  would be  to detect the `true' contamination level, that is, the minimal level of trimming for which the postulated
model would not be rejected. In this way we could, also, detect possible contaminations in the generating mechanism. 
To address this objective, we resort to the following result obtained in greater generality in  \cite{ratesConv}. 
\begin{theorem} If $\alpha\in(0,1)$ and $\nu>1$, then 
\begin{equation}
\label{dk_recorte_conv}
d_K(F,R_\alpha(F_n))=o_P\Big(\frac{(\log{n})^\nu}{n}\Big).
\end{equation}
\end{theorem}
Therefore, if $F =(1-\alpha_0)F_0+\alpha_0 G_0$ and we test for $\alpha>\alpha_0$, as $n \to \infty,$ trimming $\alpha$ from $F_n$ 
will  eliminate the part of the sample coming from $G_0$,  but also will affect the part 
of the sample coming from $F_0$.   This fact and Proposition \ref{consistency} lead to the following statement.
\begin{prop}
Let $\rho_n=O(n^{-1/2})$ and $\rho_n^{-1}=O(n^{1/2})$, and $\alpha>\alpha_0$. Then: 
\begin{equation}\label{eq_overf}
 \frac{d_K(F_0,R_\alpha(F_n))}{\rho_n}\rightarrow 
\begin{cases}
\infty \mbox{ almost surely, \hspace{.5cm} if } d_K(F_0,R_\alpha(F))>0\\
0 \mbox{ in probability,  \hspace{.7cm} if } F_0\in R_{\alpha_{0}}(F).
\end{cases}
\end{equation}
\end{prop}

\begin{figure}[htb]
\begin{center}
\includegraphics[scale=0.4]{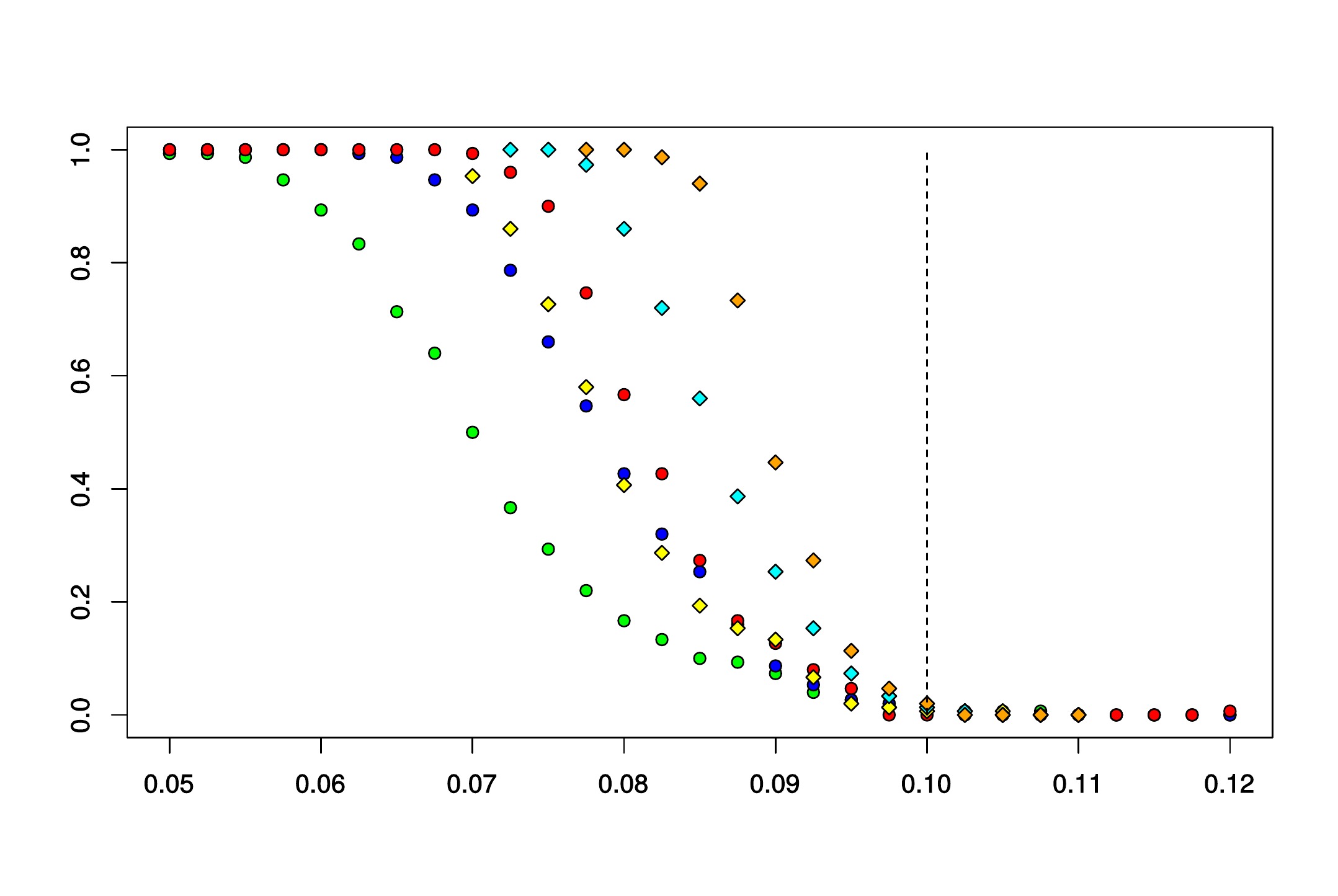}
\caption{Round green (blue, red) dots represent the frequency of rejection (y label) for 150 independent samples 
of a generating mechanism $F_1\sim 0.9N(0,1)+0.1N(3,1)$ for sample sizes 2000 (4000, 6000) and a model $F_0\sim N(0,1)$, 
as we vary the trimming level $\alpha$ (x label). Diamond yellow (cyan, orange) dots represent the rejection frequency for a generator $F_2\sim 0.9N(0,1)+0.1N(0,0.1^2)$ for sample sizes 12500 (25000, 50000). The black dashed line represents the true  contamination level  which is 0.1, since $F_0\in R_{0.1}(F_1)$ and $F_0\in R_{0.1}(F_2)$. The error probabilities are fixed to $\epsilon_1 = \epsilon_2 = 0.05$.}
\label{fig_trimm_est}
\end{center}
\end{figure}

This means that, for big enough samples, our  testing procedure will be able to detect the \emph{overtrimming} boundary, that is, the trimming level beyond which
the trimmed sample is closer to the model than true random samples from that model. In Figure \ref{fig_trimm_est} we 
are able to appreciate this behaviour (see the caption for details). The frequency of rejecting the null, for both models, 
after trimming 0.11 or more is almost zero, the theoretical contamination being 0.1. We see that around 0.1  
the models start dropping abruptly the rejection level, but  that for the model contaminated with a $N(3,1)$ we need much less points to attain the expected
behaviour than we need for the model contaminated with a $N(0,0.1^2)$. In other words,
the presence of a meaningful outlier contamination, even when trimming is allowed, disturbs more heavily the Kolmogorov distance than the presence of  equally meaningful 
inlier contamination. 
In any case, these  results suggest that it may be possible to find an estimator for the `true' contamination level. We elaborate a little bit more about this in the next section.

%
%\begin{figure}[ht]
%\centering
%\def\svgwidth{11cm}
%%\includegraphics{females_height.eps}
%\input{trimming_estimation_v2.eps_tex}
%\caption{Round green (blue, red) dots represent the frequency of rejection (y label) for 100 independent samples of a generating mechanism $F_1\sim 0.885N(0,1)+0.115N(3,1)$ for sample sizes 2000 (4000, 6000) and a model $F_0\sim N(0,1)$, as we vary the trimming level $\alpha$ (x label). Diamond yellow (cyan, orange) dots represent the rejection frequency for a generator $F_2\sim 0.828N(0,1)+0.172N(0,4)$ for sample sizes 10000 (15000, 20000). The black dashed line represents the true  contamination level  which is 0.1, since $F_0\in R_{0.1}(F_1)$ and $F_0\in R_{0.1}(F_2)$. The error probabilities are fixed to $\epsilon_2 = 0.05$ and $\epsilon_1 = 0.05/(0.999e^2)\approx 0.0068$}.
%\label{fig_trimm_est}
%\end{figure}

\section{A central limit theorem with applications}\label{Comparison} 
We divide this section in two subsections, respectively devoted to the presentation of results and to some of 
their applications. In particular, we stress on the extension of some of the applications that 
 \cite{modelAdec} and \cite{parMult} explored just on  multinomial models.

\subsection{A central limit result}\label{CLT}

What follows is our main theoretical result which describes the asymptotic behaviour of the normalized difference 
between the empirical estimator and the theoretical trimmed Kolmogorov distance under some regularity assumptions.
We recall from Section 2 that $d_K(F_0,R_\alpha(F))$ can be expressed in terms of $H^{-1}:=F_0\circ F^{-1}$.
We need to introduce the following sets, with $G$, $U$, $L$ and $\tilde{h}_\alpha$ standing for the same objects as in Theorem \ref{prop_disting_h} in Section 2,
\begin{eqnarray}\label{T1}T_{1}&=&\Big\{t \in [0,1]:\, G(t)=\|\tilde{h}_\alpha-G\|,{\textstyle  \frac 1 2 (U(t)+L(t))\geq 0}\Big\},\\
\label{T2}T_{2}&=&\Big\{t \in [0,1]:\, {\textstyle -\frac{\alpha}{1-\alpha}-G(t)= \|\tilde{h}_\alpha-G\|,  \frac 1 2(U(t)+L(t))\leq \frac{-\alpha}{1-\alpha}}\Big\},\\
\label{T3}T_{3}&=&\Big\{(s,t):\, 0\leq s\leq t\leq 1,\, {\textstyle \frac 1 2 (G(t)-G(s))}= \|\tilde{h}_\alpha-G\|, {\textstyle \frac 1 2 (G(t)+G(s))\in [-\frac \alpha {1-\alpha},0]}\Big\}.
\end{eqnarray}
A look at Theorem 4.1 in \cite{Hristo} shows that $T_1\cup T_2 \cup T_3\ne\emptyset$ provided $H^{-1}$ is continuous. We further denote
${T}^*_{1}=\{ t\in T_1:\, \frac 1 2 (U(t)+L(t))=0\}$, ${T}^*_{2}=\{ t\in T_2:\, \frac 1 2 (U(t)+L(t))=-\frac \alpha{1-\alpha}\}$ and 
${T}^*_3=\{(s,t)\in T_3:\,  \frac 1 2 (G(t)+G(s))\in \{-\frac \alpha{1-\alpha},0 \} \}$. To avoid pathological examples we will assume that
\begin{equation}\label{condiciontecnica}
{T}^*_{1}=\emptyset,\quad {T}^*_{2}=\emptyset, \quad {T}^*_3=\emptyset.
\end{equation}

Our last  regularity assumptions concern $H$, the d.f. of the random variable $F_0(X)$, where $X\sim F$. They allow 
the use of the strong approximation of the quantile process in the proof of the  theorem (developed in the Appendix). We assume that $H$ has a density, $h$
supported in $[a,b]$ (note that, necessarily, $[a,b]\subset [0,1]$) and either {one of}
\begin{equation}\label{regularityH}
h \mbox{ is positive and continuous on $[a,b]$,}
\end{equation}
\begin{equation}\label{regularityHa}
h \mbox{ is positive and continuous on $(a,b)$; for some $\varepsilon>0$, $T_1, T_2\subset [\varepsilon,1-\varepsilon]$, $T_3\subset [\varepsilon,1-\varepsilon]^2$.}
\end{equation}

\begin{theorem}
\label{tma_tcl}
Assume that $F_0$ and $F$ are continuous d.f.'s, that $F$ is strictly increasing and that the d.f. $H$ associated to $H^{-1}=F_0\circ F^{-1}$ satisfies (\ref{condiciontecnica}) and either
(\ref{regularityH}) or (\ref{regularityHa}). Then, 
\begin{eqnarray*}
\lefteqn{\sqrt{n}\left(d_K(F_0,R_\alpha(F_n))-d_K(F_0,R_\alpha(F))\right)}\hspace*{3.5cm}\\
&& \underset w \to {\textstyle \frac 1 {1-\alpha}} \max\Big(\max_{t\in T_1}B(t),\max_{t\in T_2}(-B(t)),\max_{(s,t)\in T_3}{\textstyle \frac 1 2({B}(t)-{B}(s))}\Big),
\end{eqnarray*}
where $B$ is a Brownian bridge on $[0,1]$.
\end{theorem}

 The limit distribution in  this result corresponds to the supremum of a Gaussian process. In fact, the index set for this process is often rather simple, consisting
of only one or two points as we show in our next example.

\begin{ex}\label{EjemploNormal2} {\em Trimmed Kolmogorov distances in the Gaussian model (cont.)} 
We revisit the cases studied in Example \ref{EjemploNormal1}. Recall that $F_0=\Phi$, $F=\Phi((\cdot-\mu)/\sigma)$ and
$H^{-1}(t)=\Phi(\mu + \sigma\Phi^{-1}(t))$. Hence $H(x)=\Phi\Big(\frac{\Phi^{-1}(x)-\mu}{\sigma}\Big)$, $0\leq x\leq 1$, which is
supported in $[0,1]$ and has a density which is positive and continuous on $(0,1)$. 
%Also, we have $w(t)=(H^{-1})'(t)=\sigma\frac{\varphi(\mu+\sigma\Phi^{-1}(t))}{\varphi(\Phi^{-1}(t))}$,
%where $\varphi$ denotes the standard normal density.

\smallskip
In the case $\sigma=1$ and $\mu>0$ the computations in Example \ref{EjemploNormal1} yield that $T_1=\{t_0\}$, $T_2=\emptyset$, $T_3=\emptyset$, with 
$t_0=\Phi\big(-\frac \mu 2+\frac 1 \mu \log(1-\alpha)\big)$. Applying Theorem \ref{tma_tcl} we obtain that
$$\sqrt{n}\big(d_K(R_\alpha(F_n), N(0,1)))-d_K(R_\alpha(N(\mu, 1)), N(0,1)))\underset w \to N\big(0,{\textstyle \frac {t_0 (1-t_0)}{(1-\alpha)^2}}\big).$$

When $\mu=0$ and $\sigma^2< 1$, writing $x_a= -\frac \Delta {2(1-\sigma^2)}$, $x_b=\frac \Delta {2(1-\sigma^2)}$ (with $\Delta=(8(\sigma^2-1)\log(\sigma(1-\alpha)))^{1/2}$),
$t_a=\Phi(x_a)$ and $t_b=\Phi(x_b)$ we get $T_1=\{t_a\}=\{1-t_b\}$, $T_2=\{t_b\}$, $T_3=\emptyset$ and Theorem \ref{tma_tcl}
yields
$$\sqrt{n}\big(d_K(R_\alpha(F_n), N(0,1)))-d_K(R_\alpha(N(0, \sigma^2)), N(0,1)))\underset w \to {\textstyle \frac 1 {1-\alpha}}\max\big({B(1-t_b),-B(t_b) }\big),$$
with $B$ a Brownian bridge.

\smallskip
Finally, if $\mu=0$ and  $\sigma > 1/(1-\alpha)$ then $T_1=T_2=\emptyset$, while $T_3=\{(t_a,t_b) \}$,
with $t_a=\Phi(x_a)$, $t_b=\Phi(x_b)$, $x_a= -\frac \Delta {2(\sigma^2-1)}$ and $x_b=\frac \Delta {2(\sigma^2-1)}$ and we obtain
$$\sqrt{n}\big(d_K(R_\alpha(F_n), N(0,1)))-d_K(R_\alpha(N(0, \sigma^2)), N(0,1)))\underset w \to  N\big(0,{\textstyle \frac {(1-t_b)(t_b-\frac 1 2)}{(1-\alpha)^2}}\big).$$
\hfill $\Box$
\end{ex}

The asymptotics showed in the previous example would allow to build asymptotic upper and lower confidence bounds for  the Kolmogorov distance between 
the random generator of the data and the set of $\alpha$-trimmings of the postulated normal model. In general, we would  not be able to describe the sets $T_i$ involved in the limit law, but Theorem \ref{tma_tcl} can be used to obtain conservative confidence bounds. Let $\beta\in (0,\frac 1 2)$ be given and write $Z_{\alpha}(F,F_0)$ for the limiting
random variable in Theorem \ref{tma_tcl}. Recall that for a Brownian bridge and $0\leq s\leq t\leq 1$ we have
$\mbox{Var}(B(t))=t(1-t)$ and $\mbox{Var}(\frac 1 2 (B(t)-B(s)))=\frac 1 4 (t-s)(1-(t-s))$. The $\beta$-quantile of $Z_{\alpha}(F,F_0)$ must be lower bounded by the 
$\beta$-quantile of the centered Gaussian r.v.'s $\frac{1}{1-\alpha}B(t)$, $t\in T_1$, $\frac{1}{1-\alpha}(-B(t))$, $t\in T_2$ and $\frac{1}{2(1-\alpha)}(B(t)-B(s))$, $(s,t)\in T_3$
(recall that at least one of $T_1,T_2, T_3$ must be nonempty). From the last variance computation we see that any of these centered Gaussian r.v.'s has variance at most
$\frac{1}{4(1-\alpha)^2}$, hence, a $\beta$-quantile lower bound is given by $\frac{\Phi^{-1}(\beta)}{2(1-\alpha)}=-\frac{\Phi^{-1}(1-\beta)}{2(1-\alpha)}$. Combining this with Theorem \ref{tma_tcl}
we see that
$$\liminf  P\left(\sqrt{n}\Big(d_K(F_0, R_\alpha(F_n))-d_K(F_0, R_\alpha(F))\Big)\geq {\textstyle -\frac{\Phi^{-1}(1-\beta)}{2(1-\alpha)}}\right)\geq 1-\beta.$$
Hence, 
\begin{equation}\label{upperconfidence}
d_K(F_0, R_\alpha(F_n))+\frac{\Phi^{-1}(1-\beta)}{2\sqrt{n}(1-\alpha)}
\end{equation}
is an upper confidence bound with asymptotic confidence level at least $1-\beta$ for $d_K(F_0, R_\alpha(F))$.

In order to get a simple and manageable lower bound for the Kolmogorov distance we need to pay attention to the worst cases inside the maximum of the limiting random variable in Theorem \ref{tma_tcl}. This means that we have to study the cases $T_1 = [0,a]$, $T_2 = [b,1]$, $T_3 = [a,b]$ where $0\leq a\leq b\leq 1$. We have the following inequalities
\begin{align*}
Z_{\alpha}(F_0,F) &= {\textstyle \frac 1 {1-\alpha}} \max\left(\max_{t\in [0,a]}B(t),\max_{t\in [b,1]}(-B(t)),\max_{s\in [a,b], t\in [s,b]}{\textstyle \frac 1 2({B}(t)-{B}(s))}\right)\\
&\leq {\textstyle \frac 1 {1-\alpha}} \max\left(\max_{t\in [0,a]}|B(t)|,\max_{t\in [b,1]}|-B(t)|, {\textstyle\frac{1}{2}}\left(\max_{t\in [a,b]}{B}(t)+\max_{s\in [a,b]}-{B}(s)\right)\right)\\
&\leq {\textstyle \frac 1 {1-\alpha}} \max\left(\max_{t\in [0,a]}|B(t)|,\max_{t\in [b,1]}|-B(t)|, {\max\left(\max_{t\in [a,b]}{B}(t),\max_{s\in [a,b]}-{B}(s)\right)}\right)\\
&={\textstyle \frac 1 {1-\alpha}}\max_{t\in [0,1]}|B(t)|.
\end{align*}
Now, denoting $\Psi(x)=P\left(\max_{t\in [0,1]}|B(t)|\leq x\right)$ the d.f. of Kolmogorov's distribution, we have 
$$\limsup  P\left(\sqrt{n}\Big(d_K(F_0, R_\alpha(F_n))-d_K(F_0, R_\alpha(F))\Big)\leq {\textstyle \frac{\Psi^{-1}(1-\beta)}{(1-\alpha)}}\right)\geq 1-\beta.$$
Hence, 
\begin{equation}\label{lowerconfidence}
d_K(F_0, R_\alpha(F_n))-\frac{\Psi^{-1}(1-\beta)}{\sqrt{n}(1-\alpha)}
\end{equation}
is a lower confidence bound with asymptotic confidence level at least $1-\beta$ for $d_K(F_0, R_\alpha(F))$.

In the following example we will show that the, arguably conservative, confidence bounds just obtained can be precise in practice. Of course, an efficient estimation of the sets $T_i$ could improve the precision of the coverage bands, but our simulations show that the rate of convergence can make highly unstable the estimation. In fact,  Theorem 3.1 in \cite{Alvarez-Esteban2016} addressed a simpler but similar problem involving the supremum of the difference of two independent Brownian bridges on the set where two d.f.'s attain their greatest distance. 

\begin{ex}\label{Ejemplo3} {\em Coverage rates for extreme cases.}
The bounds (\ref{upperconfidence}) and (\ref{lowerconfidence}) are
conservative. Nontheless, there are extreme cases for which the
bounds are (almost) optimal. We present several examples of such cases in Table \ref{rates_extreme}.
For different combinations of $F_0,F$ and $\alpha$, we give the observed coverage frequency of the confidence
bounds (\ref{lowerconfidence}) and (\ref{upperconfidence}). Figure
\ref{fig_extrem_cases} shows the he d.f.'s of some of these examples to get a better notion of the
functions of interest. For simplicity in all the considered cases we fix $F\sim U(0,1)$.
Then we consider instances of $F_0$
for which, approximately, $T_1\cup T_2\cup T_3$ equals $[0,1]$. For this, we fix $0\leq a\leq
b\leq 1$ and define the piecewise linear function
$${ F^{a,b}_0(t)=\begin{cases}t/(1-\alpha)+d_\alpha & t \in
[-(1-\alpha)d_\alpha,t_0]\\

a & t\in [t_0,t_1]\\

(t-(1+q)\alpha)/(1-\alpha) & t\in[t_1,(a+b)/2]\\

(t+q\alpha)/(1-\alpha) & t\in[(a+b)/2,t_2]\\

b & t\in[t_2,t_3]\\

(t-\alpha)/(1-\alpha)-d_\alpha & t\in[t_3,1+(1-\alpha)d_\alpha]\\

\end{cases}}$$
where we take $q \in(0,1)$ such that $d_\alpha= \frac{(2 +
q)\alpha}{2(1-\alpha)}<1$ and define $t_0 = (1-\alpha)(a-d_\alpha),
t_1 = (1-\alpha)a + (1+q)\alpha, t_2 = (1-\alpha)b-q\alpha, t_3 =
(1-\alpha)(b+d_\alpha)+\alpha$ (Figure 
\ref{fig_extrem_cases} depicts $F^{0.01,0.99}_0$ and $F^{1/3,2/3}_0$ for $q=0.1$ and $\alpha=0.05$) .
%Now, let us define $F^{a,b}_0(t)$ as: $t/(1-\alpha)+d_\alpha$ for
It is straightforward to check that $d_K(F_0^{a,b}, R_\alpha(F))=d_\alpha$, and that $T_1 = [0,t_0]$, $T_2 = [t_3,1]$ and
$T_3 = [t_1,t_2]$. We note that $T_1$ becomes close to $[0,a]$, $T_2$ to $[b,1]$ and
$T_3$ to $[a,b]$  as $\alpha\to 0$. 

For different extreme behaviour we take $F_0$ to be the
d.f. of a $Beta(1,\beta_0)$ distribution with $\beta_0$ such that
$f_0(1/2)=1/(1-\alpha)$  (this is possible for $\alpha < 0.06148$). We
obtain $d_K(F_0,R_\alpha(F))=P\left(Beta(1,\beta_0)\leq
1/2\right)-(1/2)/(1-\alpha)$, $T_1 = \{1/2\}$ and $T_2=T_3=\emptyset$. Figure \ref{fig_extrem_cases} includes the d.f. of $Beta(1,1.637464)$ (corresponding to $\alpha=0.05$).
Finally, another extreme case follows by fixing $d_\alpha\in (0,1)$ and
defining
\[{F^{0.5}_0(t) = \begin{cases}(1/(1-\alpha)+2d_\alpha)t & t\in[0,1/2]\\

((1-2\alpha)/(1-\alpha)-2d_\alpha)t+(\alpha/(1-\alpha)+2d_\alpha) &
t\in[1/2,1].\\

\end{cases}}\]
It is immediate that $d_K(F^{0.5}_0,R_\alpha(F))=d_\alpha$, $T_1 =
\{1/2\}$ and $T_2=T_3=\emptyset$. In Figure \ref{fig_extrem_cases} we included the case for $d_\alpha=0.1.$
\end{ex}

\begin{remark}
Notice that $F^{a,b}_0$ is not continuous in $(a+b)/2$ and is not
differentiable in $t_0,t_1,t_2$ and $t_3$, also, $F^{0.5}_0$ is not
differentiable in $1/2$. However, it is possible to modify these
functions  in such a way that from the point of view of simulation their behaviour becomes indistinguishable.  This is why we keep the simple versions that give a
better intuitive idea of what is happening.
\end{remark}

\begin{figure}[htb]
\begin{center}
\includegraphics[scale=0.4]{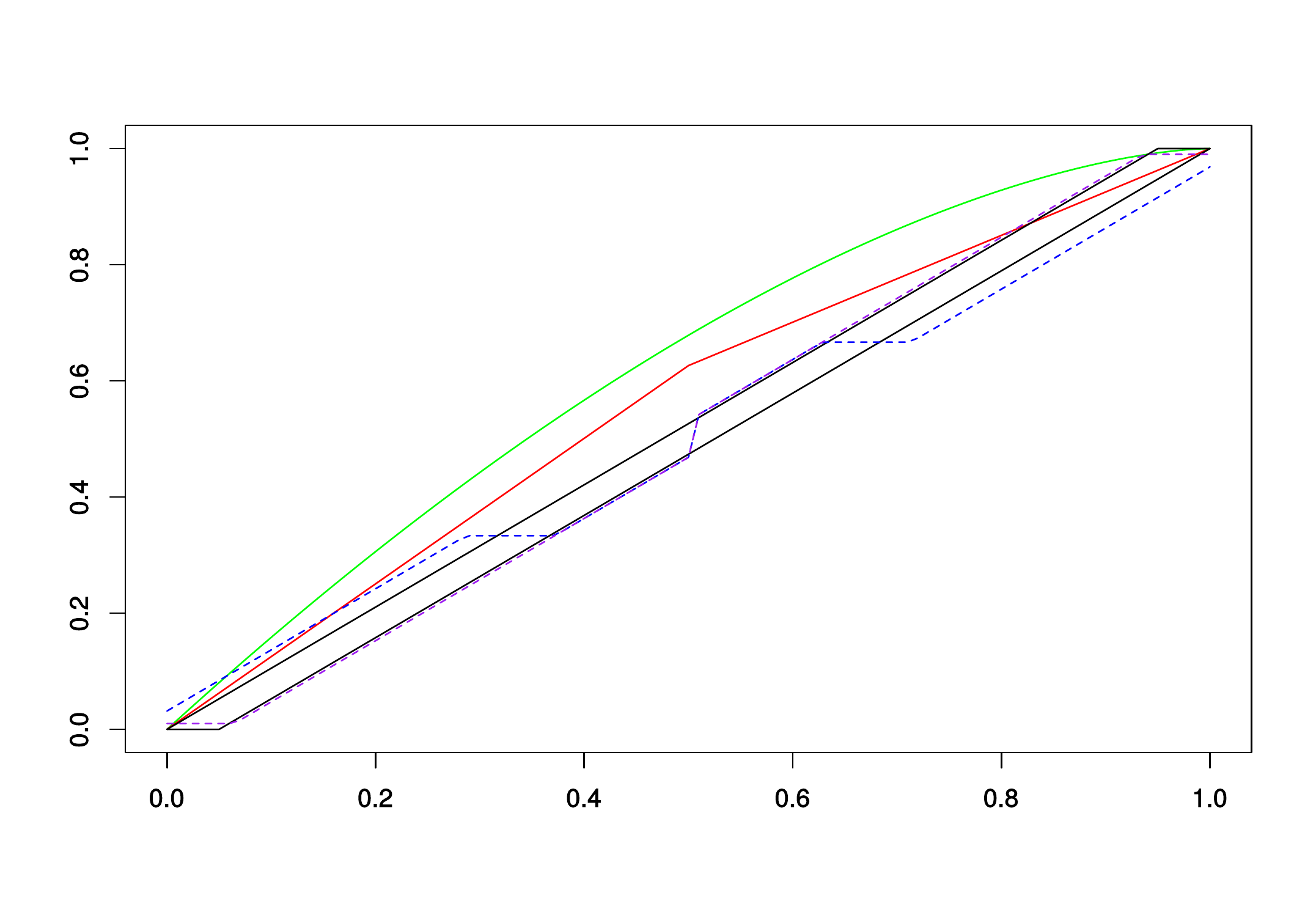}
\caption{In green $Beta(1,1637464)$; in red $F_0^{0.5}$ with
$d_\alpha=0.1$; in dashed blue $F_0^{1/3,2/3}$; in dashed purple
$F_0^{0.01,0.99}$; in black the maximum and minimum (in the usual
stochastic order sense) of the trimmings of $U(0,1)$. We fix $\alpha = 0.05$
and $q= 0.1$.}
\label{fig_extrem_cases}
\end{center}
\end{figure}


% \usepackage{multirow}
\begin{table}[ht]
\centering
\caption{For the first two examples we fix $\alpha = 0.05$ and get $\beta_0 = 1.637464$. For $F^{0.5}_0$ we take $d_\alpha=0.1$.
For all the other examples $\alpha = 0.01$ and $q = 0.01$,
where the first row indicates the values $(a,b)$ for $F^{a,b}_0$. For
each example, we generate $M=1200$ samples of size $N$ from $F \sim
U(0,1)$.}
\label{rates_extreme}
\renewcommand{\arraystretch}{1.1}
\small
\tabcolsep=0.07cm
\begin{tabular}{|c|c|c|c|c|c|c|c|c|}
\cline{2-8}
\multicolumn{1}{c|}{}&{ $Beta(1,\beta_0)$} & { $F^{1/2}_0$} & { (0.01, 0.99)} & { (0.49, 0.51)} & { (1/3, 4/3)} & { (0.01, 0.5)} & { (0.6, 0.8)} \\ \hline
\multicolumn{1}{|c|}{\multirow{2}{*}{ N = 100}}  & 0.985   & 0.993       & 0.957        & 0.955        & 0.947      & 0.949       & 0.942      \\ \cline{2-8} 
                                                                    & 0.998   & 1.000       & 1.000        & 1.000        & 1.000      & 1.000       & 1.000      \\ \hline
\multicolumn{1}{|c|}{\multirow{2}{*}{ N = 1000}} & 0.988   & 0.989       & 0.981        & 0.958        & 0.968      & 0.970       & 0.968      \\ \cline{2-8} 
                                                                    & 0.992   & 0.976       & 1.000        & 1.000        & 1.000      & 1.000       & 1.000      \\ \hline
\multicolumn{1}{|c|}{\multirow{2}{*}{ N = 5000}} & 0.993   & 0.996       & 0.998        & 0.960        & 0.973      & 0.970       & 0.958      \\ \cline{2-8} 
                                                                    & 0.980   & 0.963       & 1.000        & 1.000        & 1.000      & 1.000       & 1.000      \\ \hline
\end{tabular}
\end{table}
\subsection{Applications to credibility analysis}\label{credibilityanalysis}

{As already noted, for large enough sample sizes  a classical goodness-of-fit test would reject the null hypothesis in (\ref{test1}) and yet we could} be interested in knowing 
how well $F_0$  describes the generating mechanism behind the data.
More about this idea of resemblance, understood as similarity between generated samples and the data can be 
found in \cite{davies1995data}. In a similar spirit, a \emph{model credibility index} was introduced in \cite{modelAdec}. 
{In short, for a fixed $\delta\in (0,1)$, and a given test of fit to a model, the $\delta$-credibility index is the sample size
for which (for samples coming from the same random generator as the data) the model is rejected with probability $\delta$.
In the setting of the testing problem (\ref{test1}) with rejection rule $d_K(F_0,R_\alpha(F_n))>\lambda\rho_n$, the credibility index is the sample size $N_\delta$ (note the  dependence of $\alpha$) for which}
\begin{equation}
\label{def_credindex}
P\Big(d_K(F_0,R_\alpha(F_{N_\delta}))>\lambda\rho_{N_\delta}\Big)=P\left(d_K(F_0,R_\alpha(F_{N_\delta}))>\frac{\lambda\rho}{\sqrt{N_\delta}}\right) =\delta.
\end{equation}
{Since the underlying random generator is unknown, $N_\delta$ cannot be computed.} Sub-sampling techniques were proposed in \cite{modelAdec}, considering 
 the estimator $N_{\delta, subs}$ as the sample size such that when we take $M$ subsamples of that size the rejection frequency of the null is $\delta$.
Drawbacks of this procedure include that it is accurate only when $N_{\delta}$ is small compared to the original sample size, that $N_{\delta, subs}$ 
can never be bigger than that sample size and that the procedure is computationally demanding. We will try to address these shortcomings, while still 
maintaining the nice intuitive interpretation associated to the credibility index.

We start by writing 
\begin{eqnarray}
\nonumber
\lefteqn{P\Big(d_K(F_0,R_\alpha(F_{N_\delta}))>\lambda\rho_{N_\delta}\Big)}\hspace*{0.5cm}\\
\nonumber
&=&P\Big(\sqrt{N_\delta}\Big(d_K(F_0,R_\alpha(F_{N_\delta}))-d_K(F_0,R_\alpha(F))\Big)>\sqrt{N_\delta}\Big(\lambda\rho_{N_\delta}-d_K(F_0,R_\alpha(F))\Big)\Big)\\
\label{eq_n_star}
&=&P\left(\sqrt{N_\delta}\Big(d_K(F_0,R_\alpha(F_{N_\delta}))-d_K(F_0,R_\alpha(F))\Big)>\lambda\rho - \sqrt{N_\delta} d_K(F_0,R_\alpha(F))\right).
\end{eqnarray}
Using Theorem \ref{tma_tcl}, 
asymptotically 
we can look for 
\begin{align}
\label{p_cred_asym_2}
P\Big(\max\big(\max_{t\in T_1}\tilde{B}(t),\max_{t\in T_2}-\tilde{B}(t),{\textstyle \frac 1 2}\max_{(s,t)\in T_3}(\tilde{B}(t)-\tilde{B}(s))\big)
>\lambda\rho-\sqrt{N_\delta}d_K(F_0,R_\alpha(F))\big)\Big),
\end{align}
where we keep using $N_\delta$ for our approximation in the asymptotic regime of the original $N_\delta$.

Next, we define a lower and an upper estimate for $N_\delta$ using the probability bounds for $Z_{\alpha}(F_0,F)$ in Subsection \ref{CLT}. Thus, we define $L_\delta$ from 
$$P\left(\frac{1}{1-\alpha}\max_{t\in[0,1]}|B(t)|>\lambda\rho - \sqrt{L_\delta}d_K(F_0,R_\alpha(F))\right)=\delta$$ 
and, similarly, $U_\delta$ from
$$P\left(N\left(0,\frac{1}{4(1-\alpha)^2}\right)>\lambda\rho - \sqrt{U_\delta}d_K(F_0,R_\alpha(F))\right)=\delta.$$
{Equivalently,
\begin{equation}\label{credindex}
L_\delta=\Big(\frac{\lambda\rho - \Psi^{-1}(\delta)/(1-\alpha)}{d_K(F_0,R_\alpha(F))}\Big)^2,\, U_\delta=\Big(\frac{\lambda\rho - \Phi^{-1}(\delta)/(2(1-\alpha))}{d_K(F_0,R_\alpha(F))}\Big)^2.
\end{equation}
and it follows easily that $N_\delta\in[L_\delta, U_\delta]$. 
We see also that the empirical estimators $L_{\delta,n}$ and $U_{\delta, n}$, built replacing $F$ by $F_n$ in (\ref{credindex}), are consistent estimators of $L_\delta$ and $U_\delta$, respectively.}

We end this section discussing on the practical use of the tK-index of fit, $\alpha^*$, introduced in (\ref{eq_emp_alpha}).
{A consistent estimator $\alpha_n^*$  would provide an intuitive measure of proximity of the 
model to the data, assessing to what extent the data can be considered a contaminated sample from the model $F_0$. Recalling the setting of the testing problem (\ref{test1}) and the 
subsequent discussion,   we would reject the null hypothesis if $d_K(F_0, R_\alpha(F_n))> \lambda\rho_n$. This suggests to consider $\alpha_n^*$ as the smallest of the solutions of the equation
\begin{equation}
\label{def_alpha_star}
\sqrt{\frac{1}{2n}\log\frac{2}{\epsilon_1}} = (1-\alpha)d_K\left(F_0,R_{\alpha}(F_n)\right)\quad\text{if}\quad \sqrt{\frac{1}{2n}\log\frac{2}{\epsilon_1}}< d_K\left(F_0,F_n\right),
\end{equation}
and $\alpha_n^* = 0$ whenever $\sqrt{\frac{1}{2n}\log\frac{2}{\epsilon_1}} \geq d_K\left(F_0,F_n\right)$. This goal is feasible by numerical methods,  allowing the use of  $\alpha_n^*$ in practice. Moreover, from (\ref{def_alpha_star}) and Proposition \ref{consistency}, $\alpha_n^*$ is almost surely consistent. The carried simulations show that $\alpha_n^*$ converges rather slowly to the theoretical value. In fact, there are connections between this estimator and those considered in the FDR setting (see \cite{Genovese}), that justify this slow convergence rate even in the DCN case. Since a lower bound for $\alpha^*$ is a main goal in FDR analysis, we will deserve  some comparisons in Section \ref{FDR}.

\section{Relations with the FDR setting}\label{FDR}

To our effects, the False Discovery Rate model essentially assumes a {\bf dominated contamination model} (DCN) like (\ref{def_cont_neigh}), $F=(1-\alpha)F_0+\alpha F'$, where $F'$ (so $F$) must be stochastically dominated by $F_0$. Recall that the  stochastic order   $F' \leq_{st}F_0$ is defined by the relation $F'(x) \geq F_0(x)$ for all $x\in \mathbb R$. The DCN assumption notably simplifies the FDR theory  (which can be based on one-sided statistics), but  the methodology developed in this paper can be useful for applications in FDR in which, as often happens, the DCN can be hardly justified. To appreciate the differences between the general framework of CN's and DCN's, it seems worthwhile to take advantage of the  analyses  in Examples \ref{EjemploUniformes} and \ref{EjemploNormal1}.

\begin{ex}\label{dominatedcase} {\em Dominated contamination neighbourhoods.} {
In the scenarios considered  in Example \ref{EjemploUniformes}, only the second case of i) presents a dominated contamination   $F=(1-\varepsilon)F_0+\varepsilon F'$, with $F'\leq_{st}F_0$. In fact, between the d.f.'s $F'$ of $U(a,b)$ laws, only those verifying $a\leq \inf\{0,b\}$ and $b\leq 1$ are stochastically dominated by $F_0$. Therefore, considering 
\begin{equation}\label{domtrim}
R^-_\alpha(F,F_0):=\{F\}\cup\{F' \in R_\alpha(F): F' \leq_{st}F_0\}, 
\end{equation}
it holds $d_K(F_0,R^-_\alpha(F,F_0))=\varepsilon-\alpha$ if $F$ is the d.f. of the $U(-\varepsilon,1)$ law and $0\leq\alpha\leq\varepsilon$; $R^-_\alpha(F,F_0)=\{F\}$  for $0\leq \alpha <1$ under ii), while under iii): $R^-_\alpha(F,F_0)=\{F\}$ for $0\leq \alpha<\varepsilon$ and $R^-_\alpha(F,F_0)=\{F_0\}$ for $\varepsilon\leq\alpha<1.$ This shows that, in presence of non-dominated  contamination, trimming under the restricting domination scheme does not necessarily improve the approximation (measured under any metric).

In the Gaussian model,  stochastic dominance $N(\mu_1,\sigma_1^2)\leq_{st}N(\mu_2,\sigma_2^2)$ is equivalent to $\mu_1\leq\mu_2$ and $\sigma_1=\sigma_2$. Thus only normal distributions with $\sigma=1$ and $\mu\leq 0$ are dominated by a $N(0,1)$ law. For fixed $\alpha$, solving the relation (\ref{normaldominated})$=0$ would give the set of normal distributions that are dominated contamination versions of the $N(0,1)$ law. Therefore the only normal law in a DCN of a normal law  is the same Gaussian. In particular, in the examples considered in Figure \ref{fig_accept_region}, only the point $(1,0)$ belongs to the DCN. Of course, non-normal distributions like the mixtures $(1-\alpha)N(0,1)+\alpha N(\mu,1)$ for any $\mu<0,$ would belong to such a DCN. 
\hfill $\Box$
}
\end{ex}

Regarding the hypothesis testing problem and Theorem \ref{Testing}, note the very different nature of the problems of interest in the FDR setup: the control of the false discovery rate through a confidence lower bound and the detection of the particular false hypotheses. Resorting to a simplified version, the problem would be described through the DCN as $F=(1-\alpha)F_0+\alpha F'$, where $F_0$ is the f.d. of the $U(0,1)$ law and $F'$ is a d.f. with support on $(0,1)$ and $F'(x)\geq x$ for every $x \in (0,1)$. The null would be $\alpha=0$, and the alternative would be $\alpha>\alpha_0.$ }Acceptance of the null hypothesis with our testing procedure, for a given $\alpha$, under the DCN setting would indicate that a lower proportion than $\alpha$ false hypotheses are compatible with our data.

In the FDR setting, estimation and confidence intervals for the contamination level are  main objectives.
In fact, there are connections between the estimator defined in (\ref{def_alpha_star}) and those considered in the FDR setting (see \cite{Genovese}), that justify the  slow convergence rate even in the DCN case. Since a lower bound for $\alpha^*$ is a main goal in FDR analysis,  some comparison is in order, but previously we will introduce a new estimate.
% (see the comparative analysis in XXXX).

It is easy to see that,  (\ref{equivalencia}) is also equivalent to $P(B)\geq (1-\alpha)P_0(B)$ for any Borel set $B\subset \mathbb R$ and to $P(B)\leq (1-\alpha)P_0(B)+\alpha$ for any such set. Moreover,  the Borel sets in $\mathbb R$ can be arbitrarily well approximated by finite unions of disjoint intervals. From these considerations, we could use of the  bound
\begin{equation}\label{Kuiperequiv}
\alpha \geq \alpha(P,P_0):=1-\inf \left\{ \frac{P(J)}{P_0(J)}, \ J \mbox{ intervals in }\mathbb R\right\},
\end{equation}
noting that $\alpha(P,P_0)$ is a semicontinuous statistical functional in the sense of \cite{Donoho}, 
allowing the obtention of nontrivial lower confidence bounds for $\alpha$. This suggests that the combination of CN with the distance of Kuiper, $d_{{Kuiper}}(P,Q):=\sup\{|P(J)-Q(J)|, J \mbox{ interval in } \mathbb R\},$ could be more natural that the $d_K$ distance. That goal deserves future work, but now we devote some attention to another, novel (Bonferroni type) lower confidence bound for $\alpha(P,P_0)$, thus for $\alpha^*$:
\begin{equation}\label{Kuiper}
\hat \alpha_k= 1-\min_{(i,j)}\frac {\beta^{-1}_{j-i,n+1-j+i}(1-\gamma/M_n)}{P_0(\left[X_{(i)},X_{(j)}\right])}.
\end{equation}
Here $1-\gamma$ is the confidence level, $\beta_{k,l}^{-1}$ denotes the quantile function of the $Beta(k,l)$ distribution, and the minimum is taken over all $M_n=n(n+3)/2$ index pairs $(i,j)$ such that $0\leq i<j\leq n+1, j-i\leq n.$

{Although not implemented here, we should mention that the bound could be refined in two ways: Replace the Bonferroni quantiles
$$\beta^{-1}_{j-i,n+1-j+i}(1-\gamma/M_n)=1- \beta^{-1}_{n+1-j+i,j-i}(\gamma/M_n)$$
with $$\beta^{-1}_{j-i,n+1-j+i}(1-\gamma_n)=1- \beta^{-1}_{n+1-j+i,j-i}(\gamma_n),$$ where $\gamma_n$ is the exact $\gamma$-quantile of the distribution of
\begin{eqnarray*}&&1-\max_{0\leq i<j\leq n+1:j-i\leq n}\beta^{-1}_{j-i,n+1-j+i}(U_{(j)}-U_{(i)}) \\ &&=\min_{0\leq i<j\leq n+1:j-i\leq n}\beta^{-1}_{n+1-j+i,j-i}(1-U_{(j)}+U_{(i)})
\end{eqnarray*}
with the order statistics $0=U_{(0)}<U_{(1)}<\dots<U_{(n)}<U_{(n+1)}=1$ of a random sample from the $U([0,1])$ distribution. Furthermore, since the small intervals are more important than the large ones, one could restrict attention to all pairs $(i,j)$ of indices $0\leq i<j\leq n+1$ such that $j-i \leq d_n,$ with $d_n=\lfloor n/2\rfloor,$ say. This means, one would consider $M_n=((2n+3)d_n-d_n^2)/2$ pairs $(i,j).$}

\begin{figure}[htb]
\begin{center}
\includegraphics[scale=0.35]{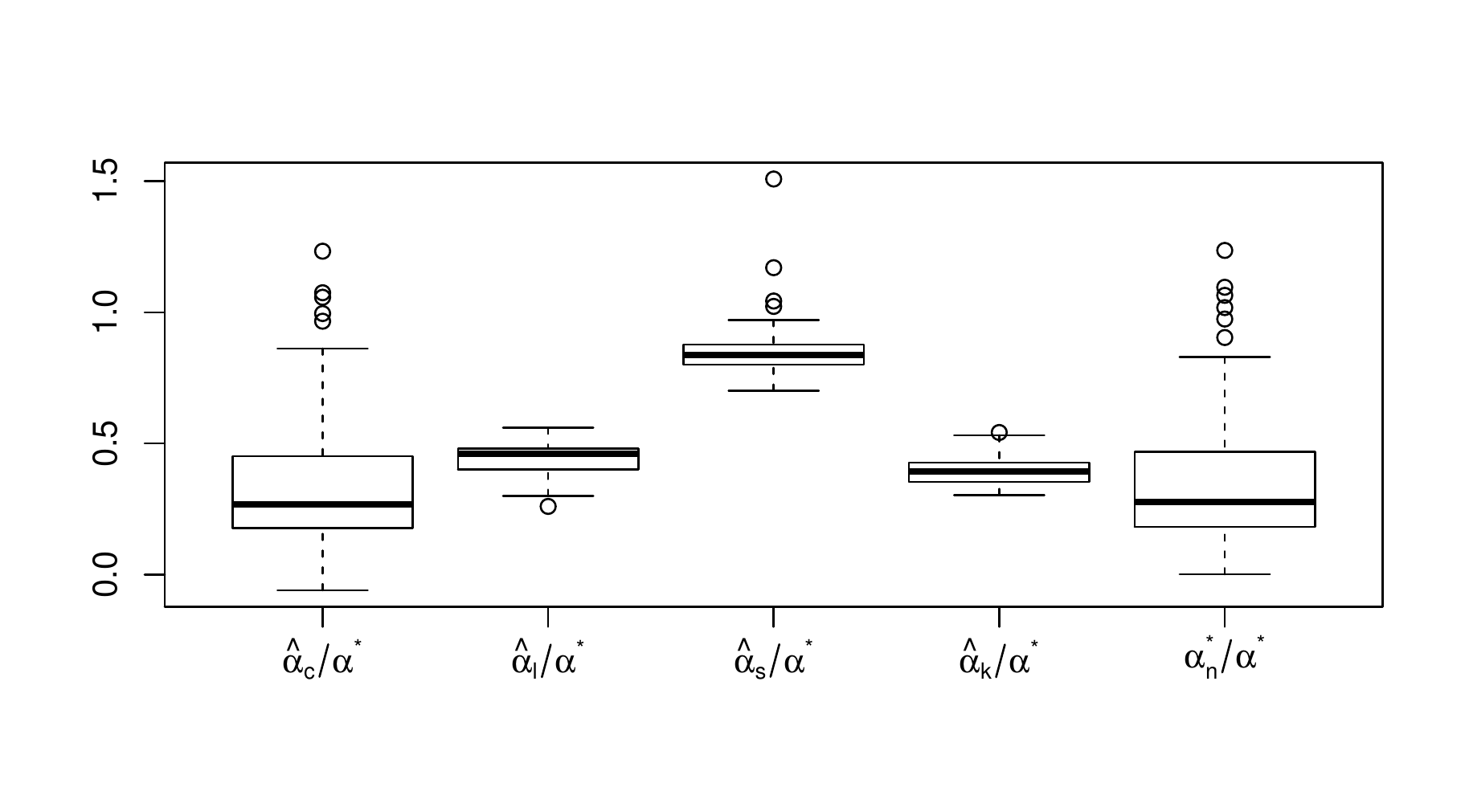}
\includegraphics[scale=0.35]{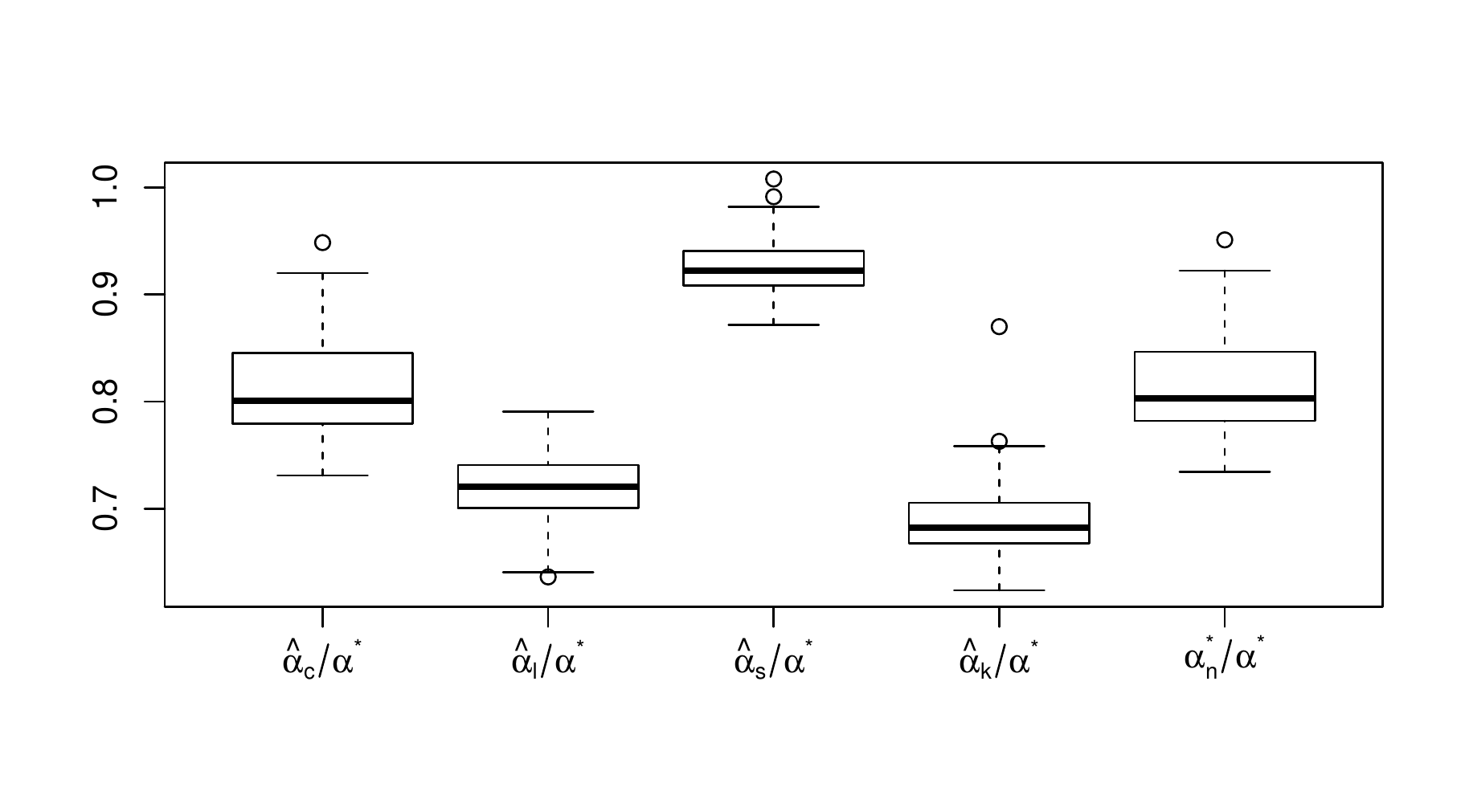}
\includegraphics[scale=0.35]{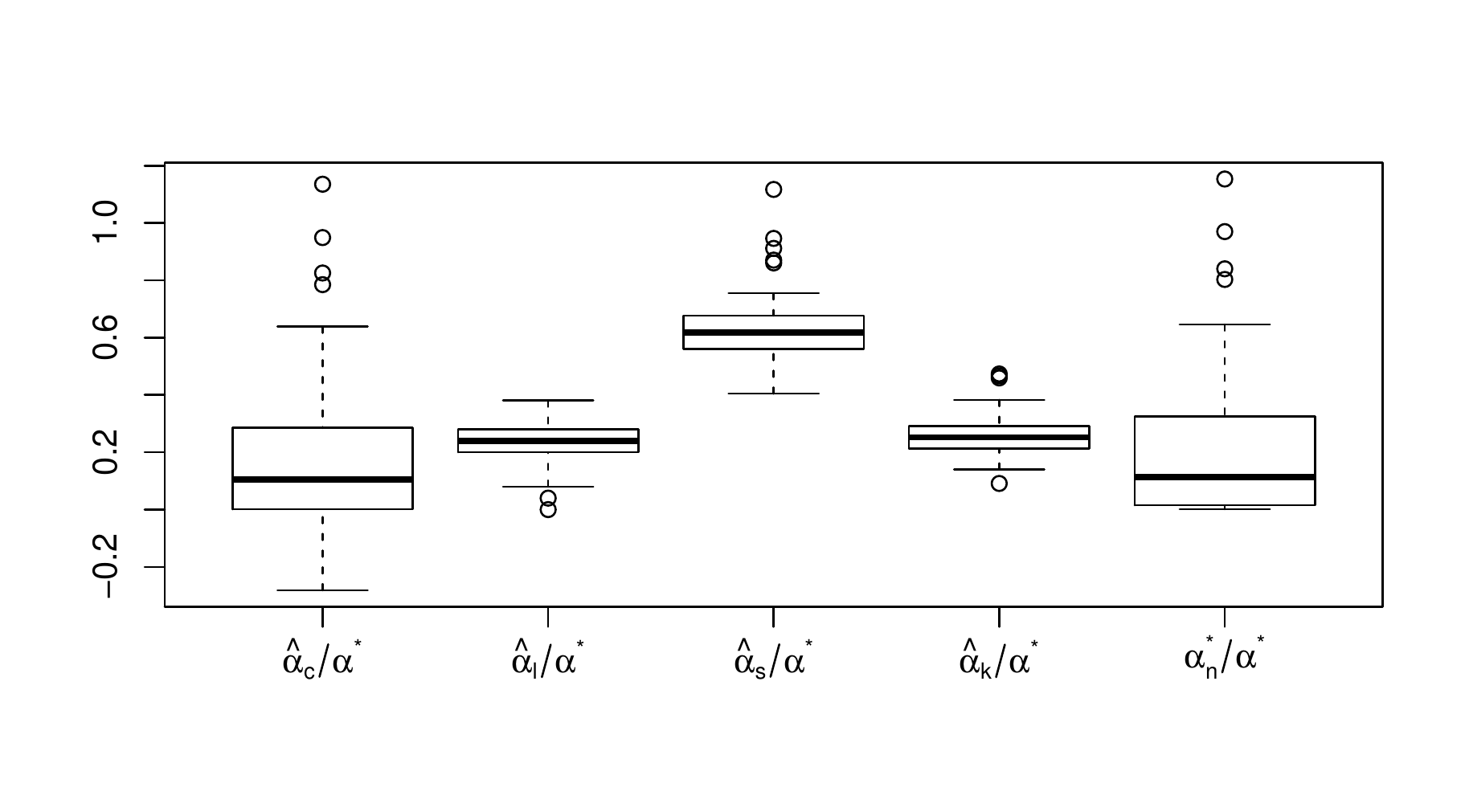}
\includegraphics[scale=0.35]{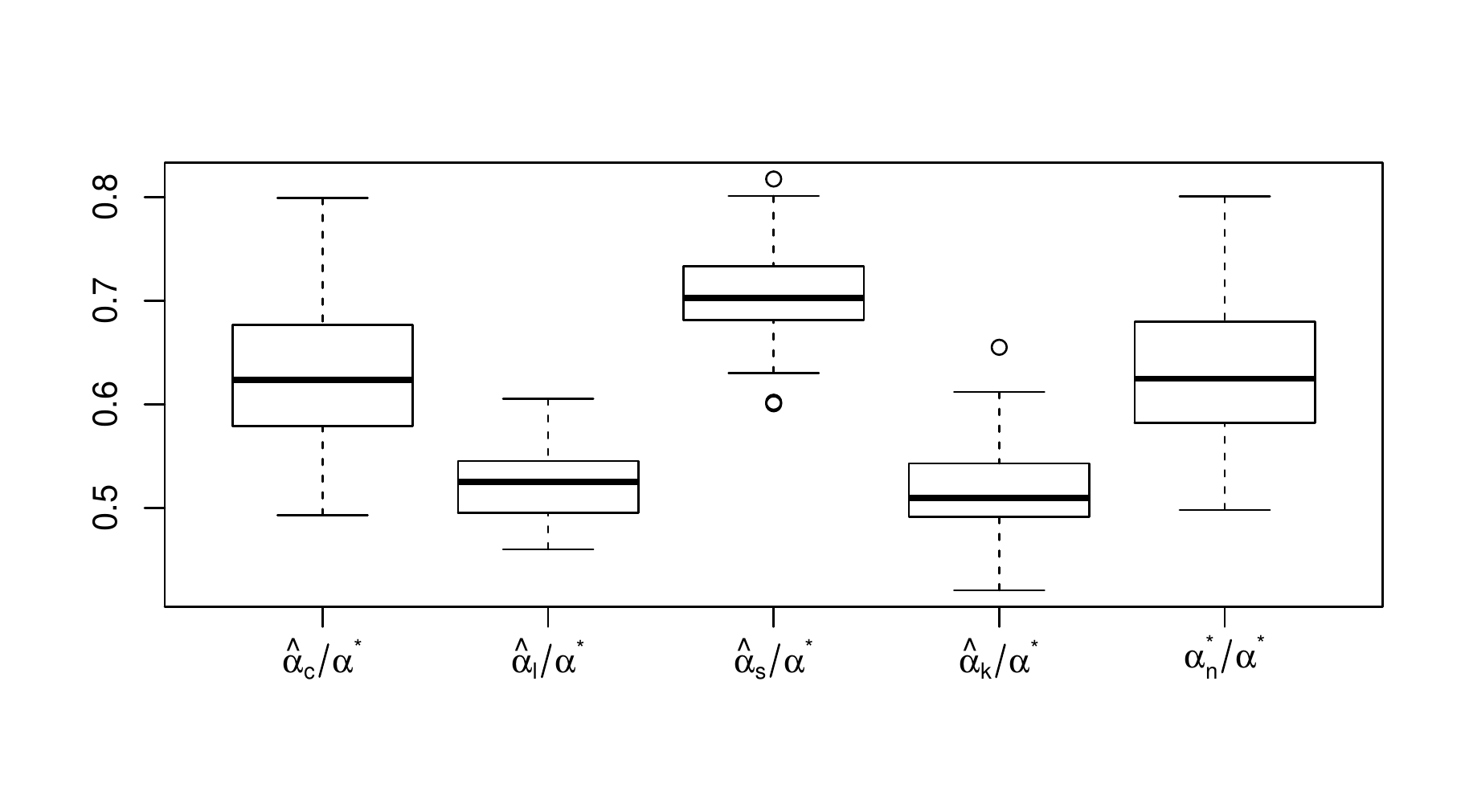}
\includegraphics[scale=0.35]{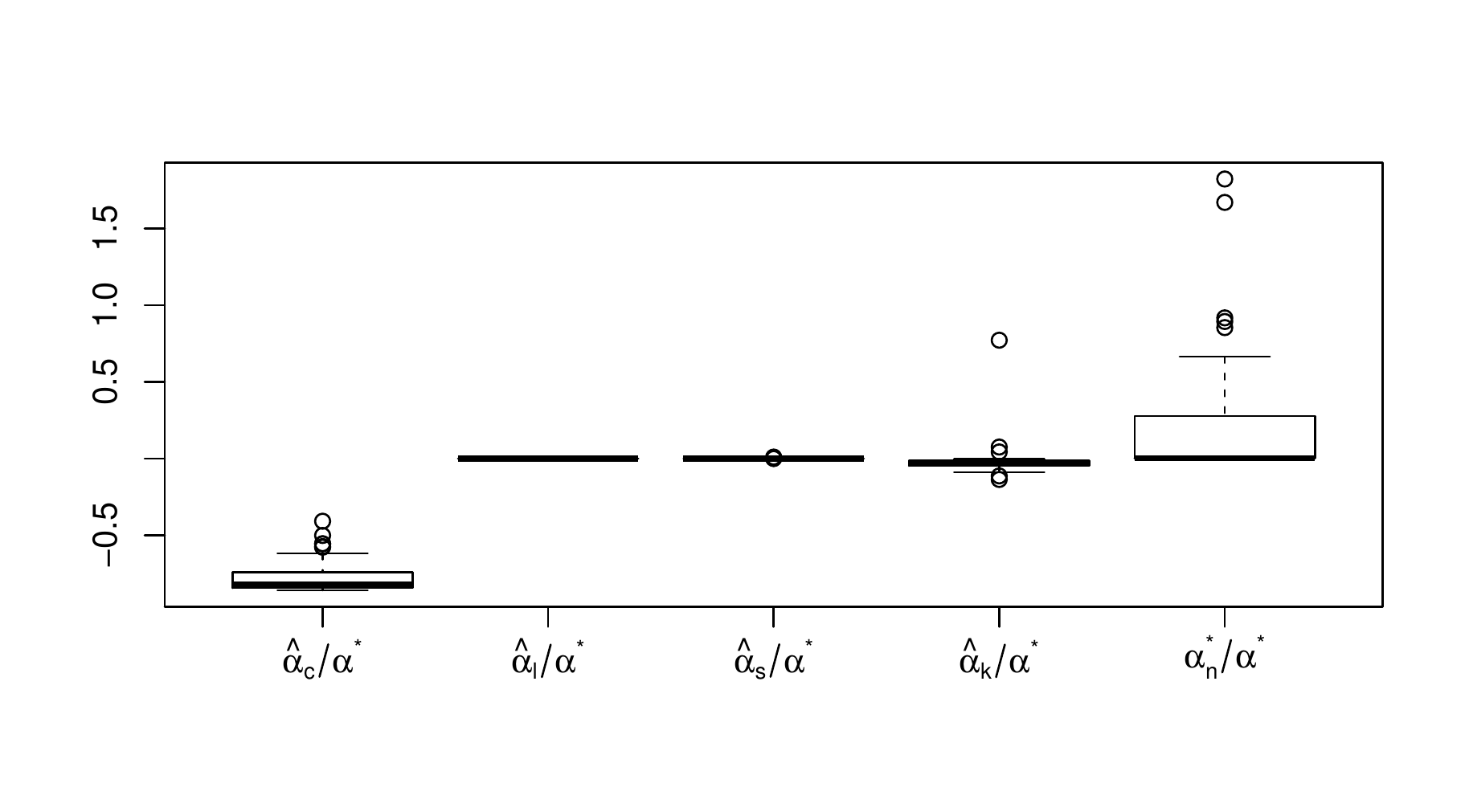}
\includegraphics[scale=0.35]{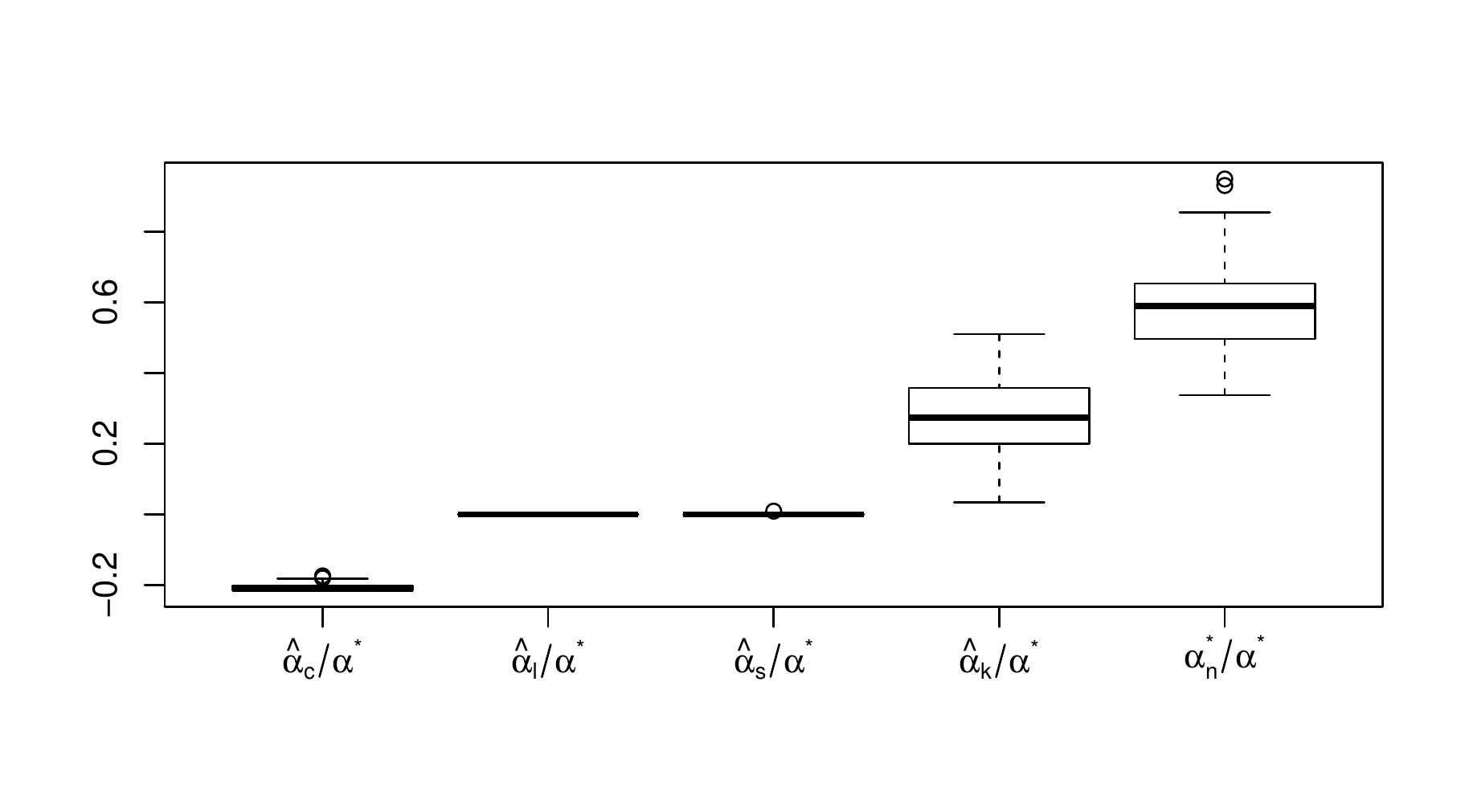}
\caption{In the graphics on the left (resp. right) column we use $\alpha^*=0.05$ (resp. $\alpha^*=0.2$). Every graphic is based on 100 samples of size 1000 obtained joining independent samples, one of size $1000(1-\alpha^*)$ of a random variable $X_i$ and other of size $1000\alpha^*$ of another $Y_i$ for respective rows $i=1,2,3$ . The estimates $\alpha_n^*, \hat\alpha_c, \hat \alpha_l, \hat \alpha_s, \hat \alpha_k$ and laws of $X_i$ and $Y_i$ are described in Example \ref{comp}}
\label{default}
\end{center}
\end{figure}

\begin{ex}\label{comp}{\em Some comparisons between estimates of $\alpha^*$.}
In Figure \ref{default}  we compare the behaviour of our estimate $\alpha_n^*$ (based on $\epsilon_1 = 0.05$) of $\alpha^*$ with some confidence lower bounds,  associated to  bounding functions, as described in \cite{Meinshausen}. We denote by $\hat\alpha_c$ to the lower bound with confidence level 0.95, i.e., $P(\hat\alpha_c\leq\alpha)\geq 0.95$, associated to the constant bounding function $\delta(t)=1$; $\hat \alpha_l$ is the one associated to the linear bounding function $\delta(t)=t$; $\hat\alpha_s$ is obtained with the standard deviation-proportional bounding function $\delta(t)=\sqrt{t(1-t)}$. The legend in the figure explains the way in which the corresponding samples have been obtained. Let $X_0$ be a random variable with a $N(0,1)$ law and recall that $\Phi$ denotes its d.f..  In the graphics of the first  row, we take $X_1=\Phi(X_0)$ and $Y_1=\Phi(X_0+4)$; in those of the second row, $X_2=\Phi(X_0), Y_2=\Phi\left(3X_0+4\right)$.  In the third row, we consider $X_3$ with a $U(0,1)$ law and $Y_3$ with a $Beta(5,1)$ law.

The first row is a very favourable case for the procedures shown in Section 4 in \cite{Meinshausen}. The second is a perturbation of that, allowing greater dispersion on the contamination, thus breaking the domination. In the lower row, we present a case where the procedures described in \cite{Meinshausen} do not give meaningful bounds, while our procedure gives sensible results. We note that \cite{Meinshausen} seemingly do not use the DCM assumption but, as it is apparent from the pictures in Figure \ref{default}, in fact their proposals are not meaningful when that condition fails.
We may conclude that our estimate is competitive when we are under the assumptions of \cite{Meinshausen}, but also works when these assumptions fail.  
\hfill $\Box$
\end{ex}

\section{Simulations and a real data example}\label{Simulations}

\subsection{A toy example}

{Let us explore  the practical use of our tools  to evaluate the 
quality of a given model on the basis of the sample. We fix the model $F_0\sim N(0,1)$
and consider three different large samples ($n=20000$) simulated from three different distributions. The first sample comes 
from $F_1\sim Logistic(0, \sqrt{3}/\pi)$,  with the same mean (=0) and variance (=1)  as   $F_0$.  We notice that this
model distribution has been reported in \cite{modelAdec} as generating datasets very close to `normality'. 
The other two samples come from contaminated normal distributions $F_2\sim 0.867 N(0, 1) + 0.133 N(0, 4^2)$ and  $F_3\sim  0.9 N(0, 1) + 0.1 N(3, 1)$ (we will refer to these samples as contaminated by `inliers' and `outliers',
respectively).

Our first step is to assess whether these samples can be assumed as coming from a contamination of level at most $\alpha=0.05$ of the model $F_0$.
We fix $\epsilon_2 = 0.05$ and $\epsilon_1 = 0.05/ (0.999e^2)$. We note (recall the discussion about the testing problem (\ref{test1}))
that this choice of $\epsilon_1$ is very close to the minimal admissible value for the validity of (\ref{rhon_explicitb}) and (\ref{rhon_explicit}), that is,
we are taking a very conservative approach, rejecting the null only if we have very strong evidence against it. From  (\ref{rhon_explicitb}) and (\ref{rhon_explicit})
we see that this amounts to fixing $\lambda=0.9997$ and $\rho_{2\times 10^4}=0.012555$, the null being rejected if 
$d_K(F_0, R_{0.05}(F_{i, 2\times10^4}))>\lambda\rho_{2\times 10^4}=0.012552$.
The first column in Table \ref{table_toy} reports the observed values of $d_K(F_0, R_{0.05}(F_{i, 2\times10^4}))$, $i=1,2,3$.
Despite the very conservative approach taken, the null is rejected for the three samples, that is, 
we should not consider them as ($0.05$) contaminated 
samples from our model $F_0$.}

\begin{table}[ht]
\renewcommand{\arraystretch}{1.5}
\centering
\caption{For  $F_0\sim N(0,1)$, $F_1\sim Logistic(0, \sqrt{3}/\pi)$, $F_2\sim 0.867 N(0, 1) + 0.133 N(0, 4^2)$ 
and $F_3\sim  0.9 N(0, 1) + 0.1 N(3, 1)$, the table shows the results obtained from samples of size $n=20000$. 
We denote $d_{K,n} = d_K\left(F_0,R_{0.05}(F_{i,n})\right)$,  $d_{K,95\%}$ are the 95\% lower (top) and upper (bottom) confidence bounds for $d_K\left(F_0,R_{0.05}(F_{i})\right)$.}
\label{table_toy}
\small
\begin{tabular}{c|c|c|c|c|c|c|c|}
\cline{2-8}
                                             & { $d_{K,n}$}                  & { $d_{K,95\%}$} & { $N_{0.5,indep}$}        & { $L_{0.5,n}$}       & { $U_{0.5,n}$}        & { $N_{0.5,subs}$}         & { $\alpha_n^*$}             \\ \hline
\multicolumn{1}{|c|}{\multirow{2}{*}{ $F_1$}} & \multirow{2}{*}{0.0140} & 0.0000         & \multirow{2}{*}{12370} & \multirow{2}{*}{4170} & \multirow{2}{*}{16079} & \multirow{2}{*}{15670} & \multirow{2}{*}{0.054} \\ \cline{3-3} 
\multicolumn{1}{|c|}{}                       &                         & 0.0262         &                        &                       &                        &                       &                        \\ \hline
\multicolumn{1}{|c|}{\multirow{2}{*}{ $F_2$}} & \multirow{2}{*}{0.0200} & 0.0000         & \multirow{2}{*}{7610}  & \multirow{2}{*}{2045} & \multirow{2}{*}{7886}  &  \multirow{2}{*}{6840}  & \multirow{2}{*}{0.069} \\ \cline{3-3} 
\multicolumn{1}{|c|}{}                       &                         & 0.0322         &                        &                       &                               &                        &                        \\ \hline
\multicolumn{1}{|c|}{\multirow{2}{*}{ $F_3$}} & \multirow{2}{*}{0.0477} & 0.0275         & \multirow{2}{*}{1135}  & \multirow{2}{*}{359}  & \multirow{2}{*}{1386}    & \multirow{2}{*}{1020}  & \multirow{2}{*}{0.089} \\ \cline{3-3}
\multicolumn{1}{|c|}{}                       &                         & 0.0599         &                        &                       &                        &                         &                        \\ \hline
\end{tabular}
\end{table}
{Next, we try to assess the quality of the rejected model as a good description of the underlying distributions of the samples. 
The simplest approach could be to use the estimated $d_K$ distance. Looking back at the first column of Table \ref{table_toy} we see 
that $F_1$ is closer to the rejection boundary than $F_2$, and the later is closer than $F_3$. 
This estimate is complemented by the $95\%$ lower (top cell) and upper (bottom cell) confidence bounds for 
$d_K(F_0, R_{0.05}(F_{i}))$, included in the second column of the table. We see, for instance, that the generator of the first sample is (with $95\%$ confidence)
at small $d_K$ distance ($0.0262$) from an $\alpha$-contamination of the standard normal distribution.

Alternatively, we could consider credibility indices, looking for the sample sizes from 
the generators that will be suitably represented by the model plus the corresponding CN (we keep our choice of $\alpha=0.05$). 
The estimators $L_{0.5,n}$ and $U_{0.5,n}$ are reported in the fourth and fifth column of Table \ref{table_toy}, we expect the 
credibility index to be in the interval $[359,1386]$ for $F_3$, in $[2045,7886]$ for $F_2$ and in $[4170,16079]$ 
for $F_1$. Once more $F_1$ is the closest to the model in this sense, followed by $F_2$ and then $F_3$. We also see 
that, from a conservative point of view, $F_1$  can generate samples of size up to 4170 while rejecting the null less 
than 50\% of the time. Therefore, at least in this sample size range $F_0$ can be considered as a useful model for the data (allowing $5\%$ 
contamination). 

In this controlled  setup we can use our knowledge of the underlying distributions of the samples to estimate the true
credibility index, $N_{0.5}$. The index
$N_{0.5, indep}$ denotes the sample size for which 5000 independent samples of that size from the true generator, 
give a rejection frequency of  50\%. $N_{0.5, subs}$ is the subsampling approximation to the credibility index described in Section \ref{credibilityanalysis}. 
We see in Table \ref{table_toy} that the 
interval $[L_{0.5,n},U_{0.5,n}]$ in all three cases contains $N_{0.5, indep} $ and $N_{0.5,subs}$, as expected. 

A last way of comparison is given by $\alpha^*_{n}$. As before,  $F_1$ is closest to the model  ($\alpha^*_n=0.054$), then 
comes $F_2$ ($\alpha^*_n=0.069$), and last $F_3$ ($\alpha^*_n=0.089$). This suggests that the random generators of the samples are 
not too far from  the model, $F_0$. On the other hand,
$F_0\in R_{0.1}(F_2)$ and $F_0\in R_{0.1}(F_3)$ and in both cases we have $\alpha^*=0.1$. Note, in this respect, the slow convergence of $\alpha^*_n$ showed in the last column 
of Table \ref{table_toy}.

To summarize, up to some `small'  contamination (0.05), the logistic generated sample is the closest one to normality. It is closer to 
normality than samples coming from 0.1-contaminations of the normal model. Also, scale contaminations with the same mean ($F_2$), generate samples 
that `look' more normal than location contaminations, when allowing some (0.05) trimming.   
}

\begin{figure}[htb]
\begin{center}
\includegraphics[scale=0.25]{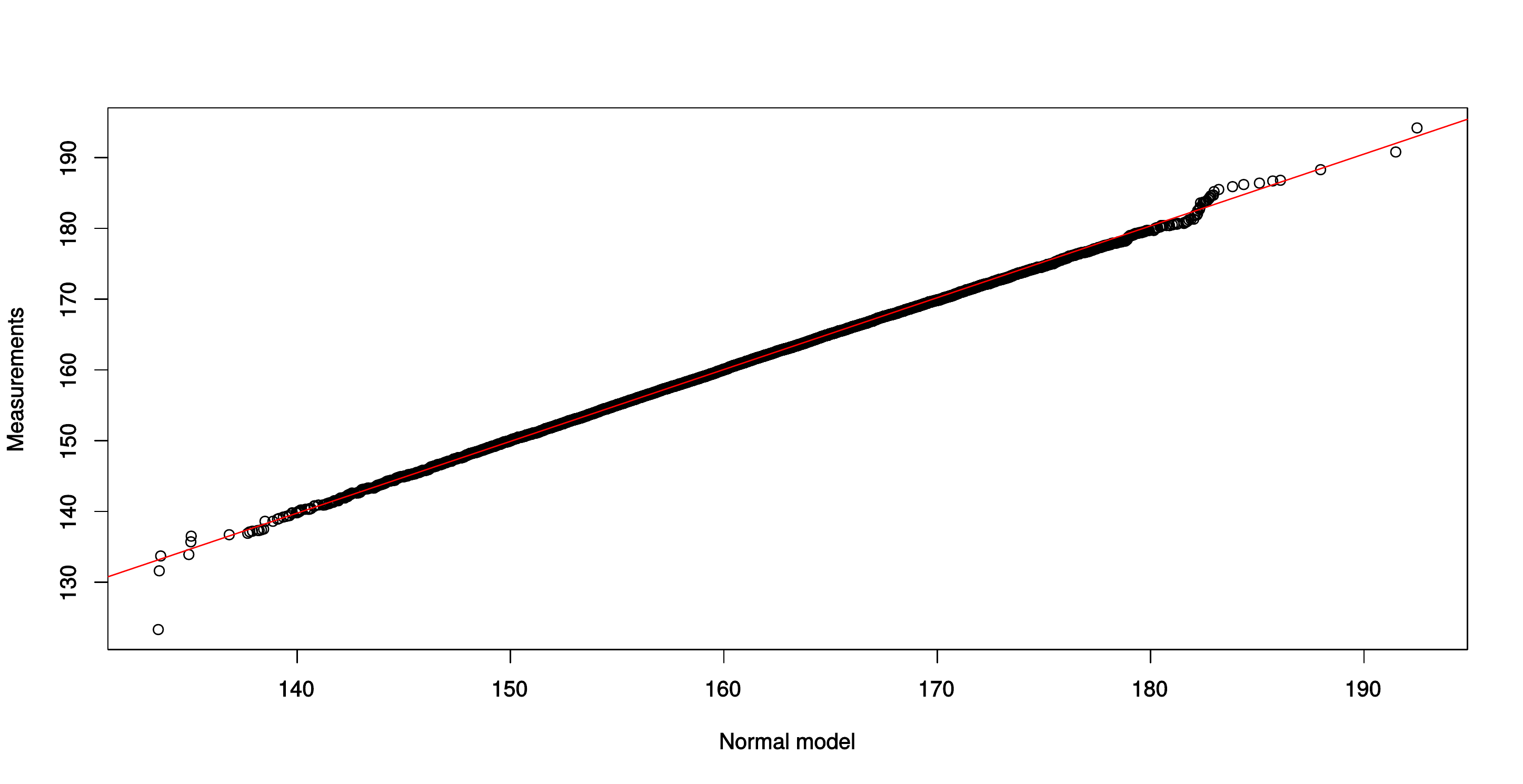}
\includegraphics[scale=0.25]{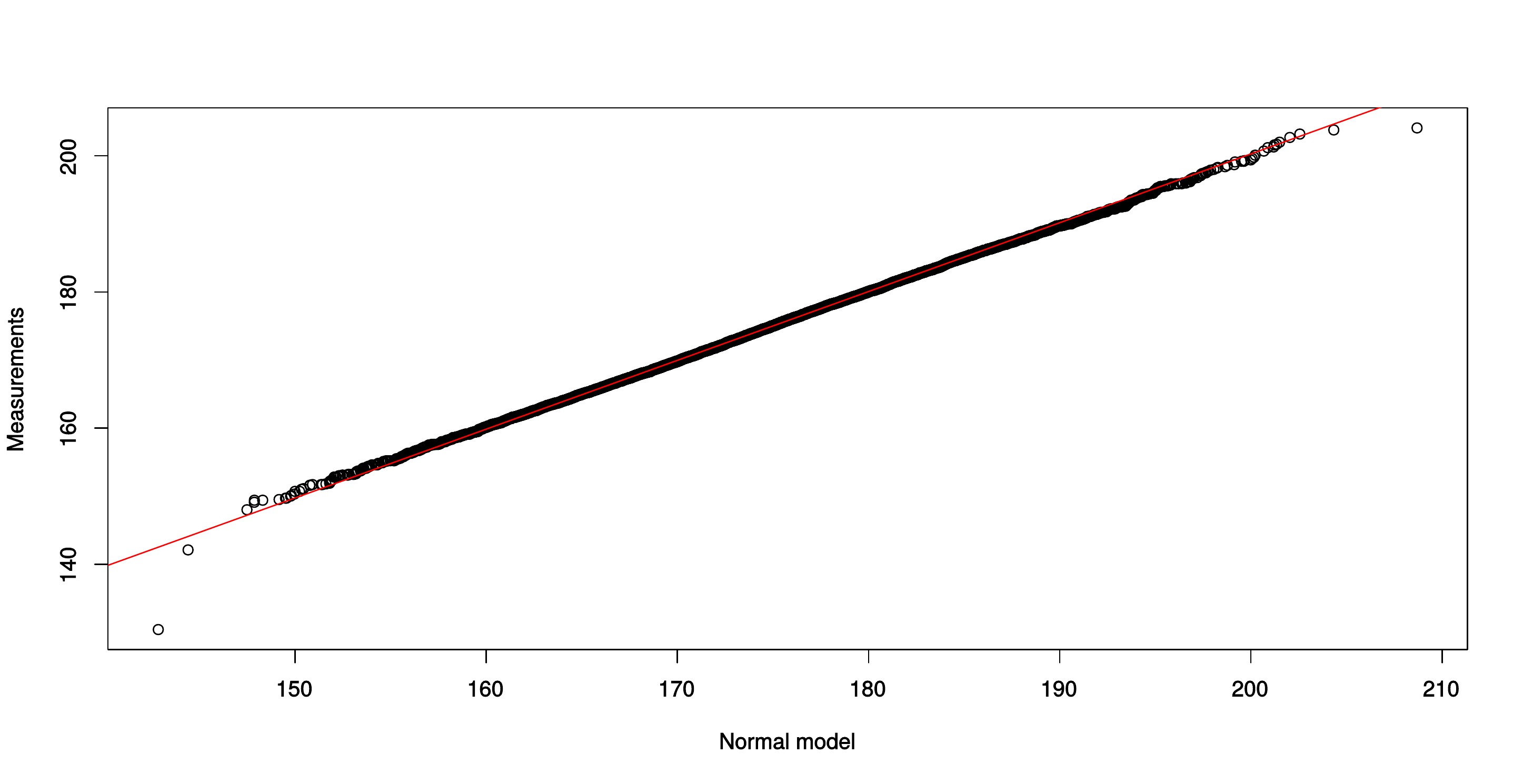}
\includegraphics[scale=0.3]{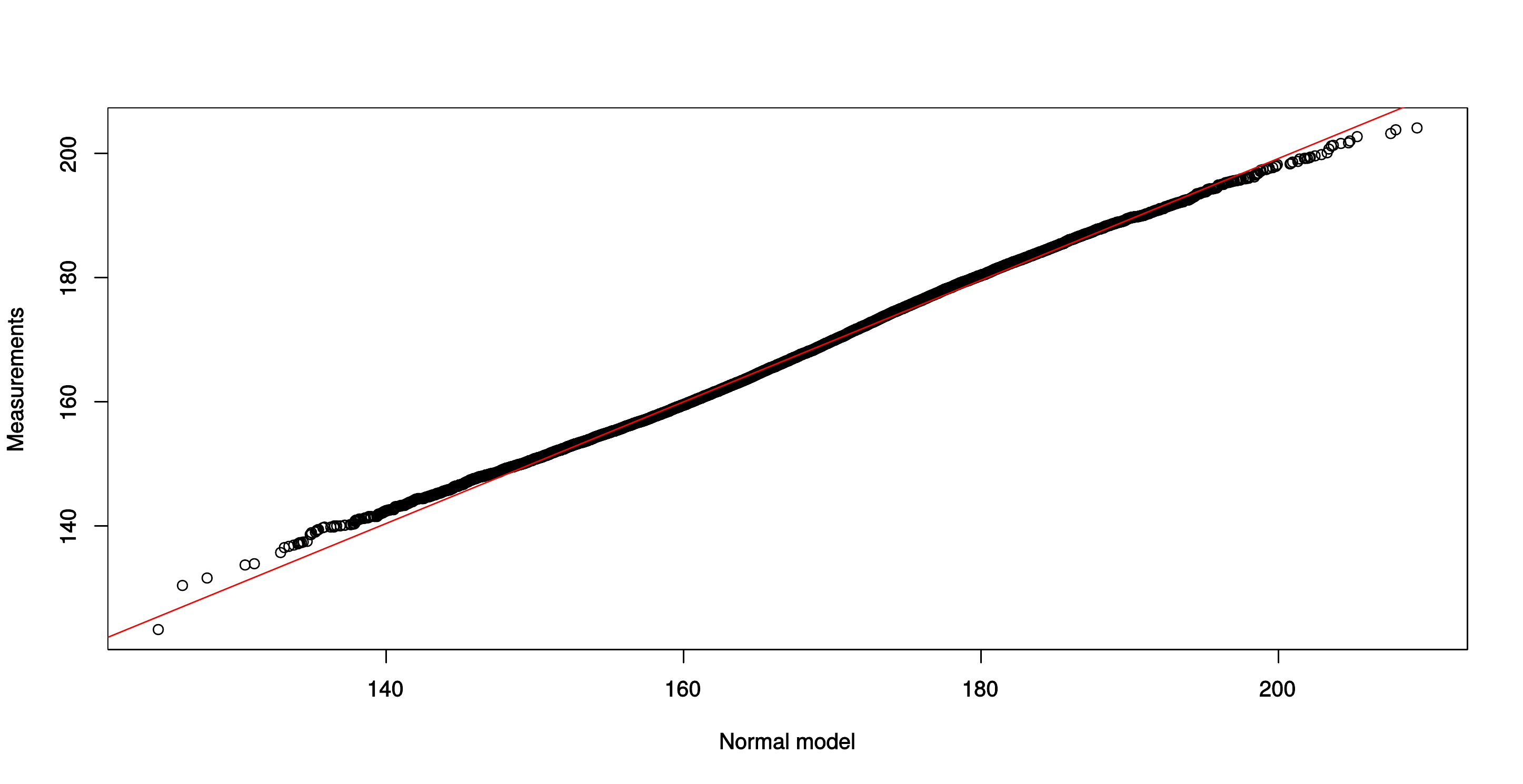}
\caption{QQ-plots of the measured heights of 15679 females (left),  14605 males (right), and of the combined joint sample (below) against a Normal distribution with the same mean and variance as the corresponding data set.}
\label{males_and_females}
\end{center}
\end{figure}

\subsection{Trying a real data example}\label{real_example}

Here we analyse the heights of 52402 individuals with ages between 2 and 84. 
The data has been obtained from NHANES (
https://www.cdc.gov/nchs/nhanes/) and consists 
of height measurements (in centimeters) of 26625 females and 25777 males. 
The dataset analysed here is available at 
\begin{center}
https://dataverse.harvard.edu/dataset.xhtml?persistentId=doi:10.7910/DVN/SHBF2G
\end{center}
We  consider  three age groups, which are related to human body development. 
The first group includes kids before puberty (ages between 2 and 10). The second group, puberty period, includes individuals  aged between 11 
and 18, with adults (over 18) making the third group. We start analysing the adult group (30284 individuals). The data consists of
height measurements on 15679 females and 14605 males.

Notice that in our analysis we will use the population estimates of the mean and variance. This is very usual in the goodness of fit setting based on procedures designed for testing simple hypothesis and, in particular in the FDR setting. There, the $U(0,1)$ law, considered as the hypothesis, arises from the integral, or $p$-value transformation, but it depends on the (unknown) true distribution. In our framework, that license is even more permisible because we are interested in getting a useful description of the data. 

We analyse first the sample by gender group. In Figure \ref{males_and_females} we see qq-plots from normal distributions with 
the same mean and variance as the female's and male's heights data. The pictures  
suggest that the normal model could provide a reasonable description of the data. 
Also, a Kolmogorov-Smirnov test of normality yields a $p$-value of $0.5385$ 
for the male group  and of $0.2997$ for the female group, thus we do not get enough evidence to reject that the data sets come from  normal
random generators.

Next we take a look at the combined data set. The previous analysis suggests the model 
$F_0\sim 0.52 N(161.0, 7.1^2)+0.48 N(174.6, 7.9^2)$. 
If, however, we perform a gender-blind analysis and take a look 
at the qq-plot in the third graphic in Figure \ref{males_and_females} for the combined sample, we may 
be tempted to say that the normal distribution is not a bad model for the data (nevertheless, 
the K-S  normality test yields a $p$-value smaller than $10^{-16}$). 
After the discussion in the previous sections, 
%and knowing that classical goodness of fit tests can be misleading for large samples, 
we could yet stick to the normal model and consider $F_0^*\sim N(167.6, 10.1^2)$
hoping a useful description of the random generator of the data.

\begin{figure}[htb]
\begin{center}
\includegraphics[scale=0.4]{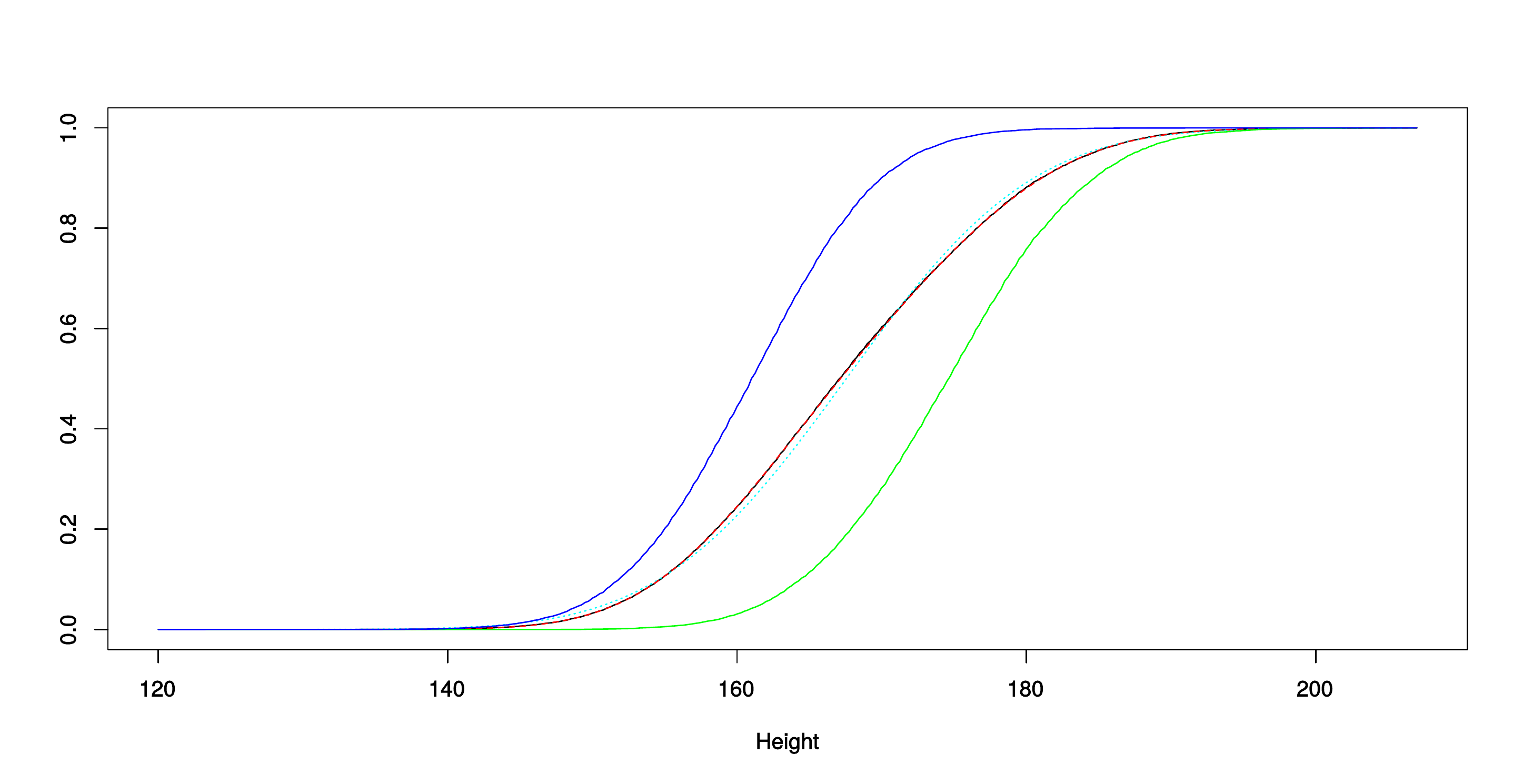}
\caption{Comparison of data and models. Solid lines correspond to the empirical d.f.'s: black for the joint samples, and blue (resp. green) for females (resp. males).  $F_0$ is represented in dashed red and   $F_0^*$  in dotted cyan.}
\label{fig_hum_alien_mod}
\end{center}
\end{figure}

Figure \ref{fig_hum_alien_mod} shows the empirical d.f. together with the models 
$F_0$ and $F_0^*$. While the gender-blind model, which is in disadvantage (since it is blinded 
to relevant information), is further away from the data than $F_0$ we may wonder how bad is $F_0^*$ as a model. 
If trimming is allowed, we would need a  $6\%$ trimming to avoid rejection of the null hypothesis, i.e., $d_K(F_0^*,R_{0.06}(F))=0$ 
would not be rejected, thus $\alpha^*_n=0.06$, and our data are compatible with a generation from $F_0^*$ with a proportion of until $6\%$ wrong data.
Actually, $F_0$ is still a better model, since  $d_K(F_0,R_{0.06}(F_n)) = 0.00231$ with $95\%$ confidence interval for 
$d_K(F_0,R_{0.06}(F))$ of $[0,0.01468]$, while $d_K(F_0^*,R_{0.06}(F_n)) = 0.01026$ with a confidence interval for $d_K(F_0^*,R_{0.06}(F))$
of $[0,0.02264]$.  We  could even look for other normal distributions inside the tolerance region, shown in green in Figure 
\ref{fig_accep_reg_aliens}, and choose one as a sensible model. Alternatively, if we  find this trimming level 
unacceptable, we may try to use smaller CN's and asses model adequacy using credibility 
analysis. 
\begin{table}[ht]
\centering
\caption{ $d_{K,\alpha, n} = d_K\left(F_0^*,R_{\alpha}(F_n)\right)$, where  $F_0^*\sim N(167.6, 10.1^2)$ and $F$ is the true generating mechanism. $N_{0.5,subs}$ is obtained taking 1000 sub-samples of the heights data. $N_{0.5, indep}$ is obtained taking 1000 independent samples from $F = F_0$.}
\label{tabel_cred_analysis}
\small
\begin{tabular}{|c|c|c|c|c|c|c|c|c|}
\hline
              {$\alpha$}                                 & {$d_{K,\alpha,n}$}        & {$d_{K,\alpha, 95\%}$} & {$N_{0.5,{\rm indep}}$}        & {$L_{0.5,n}$}           & {$U_{0.5,n}$}            & {$L_{0.5,95\%}$} & {$U_{0.5,95\%}$} & {$N_{0.5,{\rm subs}}$}         \\ \hline
\multicolumn{1}{|c|}{\multirow{2}{*}{$0.015$}} & \multirow{2}{*}{0.0184} & 0.0000               & \multirow{2}{*}{8350}  & \multirow{2}{*}{2239} & \multirow{2}{*}{8631}  & 832           & 3206           & \multirow{2}{*}{7225}  \\ \cline{3-3} \cline{7-8}
\multicolumn{1}{|c|}{}                                  &                         & 0.0302               &                        &                       &                        & $\infty$      & $\infty$       &                        \\ \hline
\multicolumn{1}{|c|}{\multirow{2}{*}{$0.035$}} & \multirow{2}{*}{0.0143} & 0.0000               & \multirow{2}{*}{17250} & \multirow{2}{*}{3888} & \multirow{2}{*}{14993} & 1143          & 4407           & \multirow{2}{*}{13280} \\ \cline{3-3} \cline{7-8}
\multicolumn{1}{|c|}{}                                  &                         & 0.0263               &                        &                       &                        & $\infty$      & $\infty$       &                        \\ \hline
\multicolumn{1}{|c|}{\multirow{2}{*}{$0.055$}} & \multirow{2}{*}{0.0110} & 0.0000               & \multirow{2}{*}{37300} & \multirow{2}{*}{6851} & \multirow{2}{*}{26417} & 1523          & 5872           & \multirow{2}{*}{25810} \\ \cline{3-3} \cline{7-8}
\multicolumn{1}{|c|}{}                                  &                         & 0.0233               &                        &                       &                        & $\infty$      & $\infty$       &                        \\ \hline
\end{tabular}
\end{table}

The output of this type of analysis is reported in Table \ref{tabel_cred_analysis}. We have fixed $\epsilon_1=0.01$ (recall 
the discussion leading to (\ref{rhon_explicit})) and considered three different trimming levels ($\alpha=0.015,
0.035$ and $0.055$, leading to optimal rejection boundaries $\lambda\rho_{30284}= 0.0098, 0.0100$ and $0.0103$,
respectively). Testing for these contamination levels results in 
rejection of the null hypothesis, since the values in the first column of Table \ref{tabel_cred_analysis} are 
above the respective rejection boundaries. But we see how the empirical trimmed Kolmogorov distance approaches the 
rejection boundary as the trimming level increases. Further informative values are provided by the intervals 
$[L_{0.5,n},U_{0.5,n}]$. With the available data, since we expect $N_\delta$ to be in $[L_{0.5,n},U_{0.5,n}]$,
we see how reasonably the gender-blind model, $F_0^*$, represents the data. With a conservative 
point of view,  after trimming only $0.055$, samples of the true generator of size $6851$ will not be rejected as coming from $F_0^*$ 
(more than $50\%$ of the time). If we take an optimistic point of view, we can say the same thing but 
for a sample size of $26417$. As in our toy example, we see that $N_{0.5, subs}\in[L_{0.5,n}, U_{0.5,n}]$, therefore our
estimated interval for the credibility index contains the estimation proposed in \cite{modelAdec}. If, on the other hand, we admit
that the data comes from $F_0$ and calculate the estimate $N_{0.5,indep}$, we see that $N_{0.5,subs}$ is far from 
$N_{0.5,indep}$ and our upper bounds $ U_{0.5,n}$ get closer to $N_{0.5,indep}$. Furthermore, we could plug-in our upper and lower confidence 
bounds (\ref{upperconfidence}) and (\ref{lowerconfidence}) into (\ref{credindex}) to get upper and lower confidence bounds for $L_{0.5}$ and $U_{0.5}$.
These are reported in the columns labeled $L_{0.5,95\%}$ and $U_{0.5,95\%}$.
We can assure with more than $95\%$ confidence that $N_{0.5}\geq 1143$ for $\alpha=0.035$ and, similarly, that $N_{0.5}\geq 1523$ for $\alpha=0.055$.

Finally, we study the normality of the data for grouping ages. Using the same mean and variance as the data, we propose 
 $F_1\sim N(116.9,18.2^2)$ for the age group under 11, $F_2\sim N(163.1,10.9^2)$ for the ages  11 and 18, 
and  $F_3 (=F_0^*)\sim N(167.6, 10.1^2)$ for ages over 18. 
%First, the empirical Kolmogorov distance and the confidence intervals give us a possible ordering: $0.01591$ and a $95\%$ confidence interval of $[0,0.02761]$; $0.02566$ and $[0,0.04402]$
The tK-index of fit allows us to compare how normal is the data in each age group. We obtain the following indices: 
$\alpha^*_{1,n}=0.3665$, $\alpha^*_{2,n} = 0.0057$ and, as before, $\alpha^*_{3,n}=0.06$. This gives a clear  
`normality' ranking. Somewhat surprisingly the data from the puberty group (ages 11 to 18) is almost normal. 
The adult group is close to normality and the children group is very far from normality. 
We emphasize that  normality  is rejected for each data set by a K-S test. To 
gain some intuition of what is really happening, we plot in Figure \ref{fig_accep_reg_aliens} the tolerance region for the normal 
family inside each respective CN for $\alpha^*_{2,n}$ and $\alpha^*_{3,n}$. The plot shows remarkably well 
how much closer to being normally distributed is the data of the teenagers compared to the adult group.

\begin{figure}[htb]
\begin{center}
\includegraphics[scale=0.4]{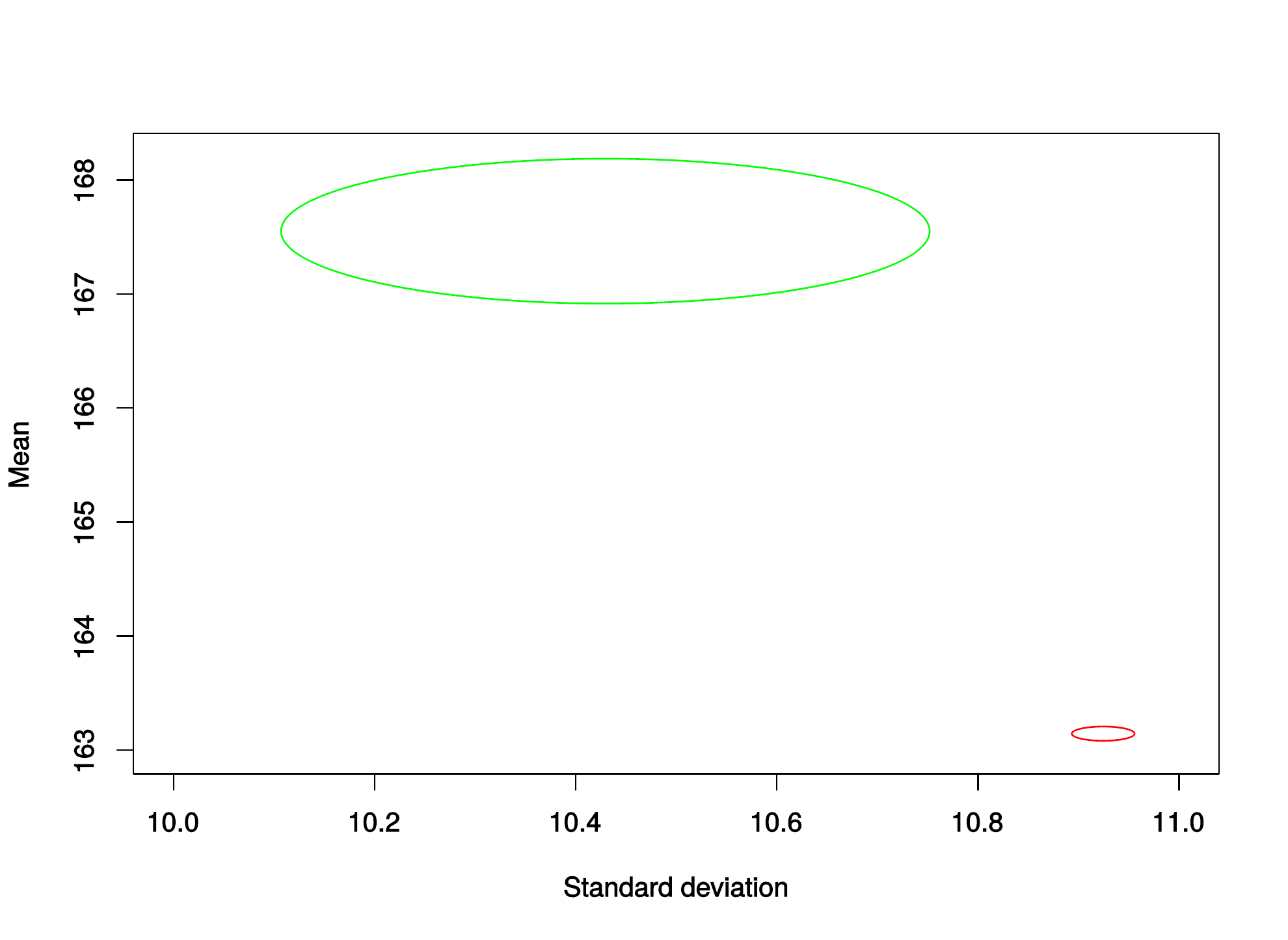}
\caption{Tolerance regions for the normal family based on $F_2$ and $F_3$ in Section \ref{real_example}. In red the tolerance region that 
is inside $\mathcal{V}_{0.0057}(N(163.1,10.9^2))$, in green the one that is inside $\mathcal{V}_{0.06}(N(167.6,10.1^2))$.}
\label{fig_accep_reg_aliens}
\end{center}
\end{figure}

\section{Conclusions}\label{Discusion}

Through the paper we showed that the Kolmogorov distance, the credibility index bounds and the tK-index of fit, provide an 
intuitive and easy to understand comparison between models. The Kolmogorov distance between a contamination model and a 
generator gives a straightforward way of comparing accepted or rejected models and, further, allows the use of the other 
two indexes in the case of rejection. The credibility index bounds provide a summary of which model is closest to the data 
and give an idea of the region in which the model agrees well with the data. The tK-index of fit provides a single
summary that can be widely used and can have attached some informative tolerance region. The procedure we have followed 
to calculate the normal family tolerance region can be more or less directly extended to other absolutely continuous distributions.
Last but not least, we have provided an efficient algorithm for computing $d_K(F_0, R_\alpha(F_n))$ which makes possible 
the implementation of all the previous procedures. 

With these tools we elaborate on the idea that rejecting a model 
does not mean that the model is useless.
Our testing procedure and limit results allow different applications of this idea. As showed in our toy example, we can use 
them to asses how some known generating mechanism produces data compatible with some fixed model when we allow some ``small" 
contamination. In this way we may obtain some useful (hopefully simpler or faster to implement) generators for some range of sample sizes. 
These tools allow also to compare different data sets, from unknown generators, to a contamination model and rank how 
well the model agrees with the data.

\section*{Appendix}

In this section we provide the proof of  Theorem \ref{tma_tcl}. The key for this result relies on 
appropriate characterizations of the best approximation of a function (in uniform norm) by monotone functions
with a box constraint, obtained in \cite{Hristo}. These characterizations allowed the obtention of a directional differentiability result (Theorem 4.3 in \cite{Hristo}) that we include below. We give a version that is simply a rephrasing  with the present notation, and involves the sets $T_i$ introduced at the beginning of Section \ref{CLT}.
As it is the case here, this kind of result typically allows to obtain efficiency and asymptotic distributional behaviour of  functionals  in the statistical setting {(see e.g. \cite{Carcamo})}.
 
\begin{theorem}\label{asymptotic1}
Assume $G,J:[0,1]\to \mathbb{R}$ are continuous functions and $r_n>0$ is a sequence of real numbers such that $r_n\to \infty$. 
Define $G_n=G+\frac J{r_n}$ and consider $U, L, \tilde h_\alpha$ as in Theorem \ref{prop_disting_h} and $\tilde h_{\alpha,n}$ built in the 
same way as $\tilde h_\alpha$ but from $G_n$. Assume further that $T_1, T_2$ and $T_3$ are as defined in (\ref{T1}), (\ref{T2}) and (\ref{T3}) and that
there is no $t\in {T}_1$ with $(L(t)+U(t))/2=0$, no $t\in {T}_2$ with $(L(t)+U(t))/2=-\alpha/(1-\alpha)$ and no $(s,t)\in {T}_3$ with $\frac 1 2 (G(t)+G(s))\in \{-\alpha/(1-\alpha),0\}$.
Then
$$r_n(\|G_n-\tilde h_{\alpha,n}\|-\|G-\tilde h_\alpha\|)\to \max\left(\max_{t\in {T}_1} J(t), \max_{t\in {T}_2} (-J(t)),\frac 1 2 \max_{(s,t)\in {T}_3} (J(t)-J(s)) \right).$$
\end{theorem}

{

\medskip 
\noindent\textit{\bf Proof of Theorem \ref{tma_tcl}.}
As in Theorem \ref{prop_disting_h} for $\Gamma$, we write $G(t)=H^{-1}(t)-\frac t{1-\alpha}$, 
$G_n(t)=H_n^{-1}(t)-\frac t{1-\alpha}$, keep the notation for $\tilde{h}_\alpha$ and write $\tilde{h}_{\alpha,n}$ for the 
corresponding object defined from $G_n$. With this notation we will show weak convergence of
$$A_n=\sqrt{n}\left(\|\tilde{h}_{\alpha,n}-G_n\|-\|\tilde{h}_\alpha-G\|\right)$$
to complete the proof. With this goal we consider the quantile process 
$$Q_n(t)=\sqrt{n}(H_n^{-1}(t)-H^{-1}(t)),\quad 0\leq t\leq 1.$$
Assumption (\ref{regularityH}) allows us to apply Theorem 18.1.1, p. 640 and Example 18.1.2, p.641, in \cite{Shorack1986} to conclude that
we can choose a version of $Q_n$ and a Brownian brigde, $B$, such that if $w=(H^{-1})'$ and $\tilde{B}=w B$ then
\begin{equation}\label{strongap}
\|Q_n-\tilde{B}\|\to 0
\end{equation}
in probability. If (\ref{regularityHa}), instead of (\ref{regularityH}), holds then we can still find versions of $Q_n$ and $\tilde{B}$
such that $\max_{\varepsilon\leq t\leq 1-\varepsilon} |Q_n(t)-\tilde{B}(t)|\to 0$ in probability. It is easy to see that (\ref{regularityHa})
implies that $$\|\tilde{h}_\alpha-G\|=\max_{\varepsilon\leq t\leq 1-\varepsilon}|\tilde{h}_\alpha(t)-G(t)|$$ and also that, in a probability one
set, eventually $$\|\tilde{h}_{\alpha,n}-G_n\|=\max_{\varepsilon\leq t\leq 1-\varepsilon}|\tilde{h}_{\alpha,n}(t)-G_n(t)|.$$ From this point we assume that 
(\ref{regularityH}) (hence, also (\ref{strongap})) holds. Our last comments, however, show that our proof can be trivially adapted to cover the case when
(\ref{regularityHa}) holds. We omit further details.

Next, we note that $w$ is a continuous function and, as a consequence, $\tilde{B}$ has, with probability one, continuous trajectories.
We note that $G_n(t)=G(t)+\frac{Q_n(t)}{\sqrt{n}}$ and introduce $\bar{G}_n(t)=G(t)+\frac{\tilde{B}(t)}{\sqrt{n}}$ and the related
functions $\bar{U}_n$, $\bar{L}_n$ and $\bar{h}_{\alpha,n}$ related to $\bar{G}_n$ as ${U}_n$, ${L}_n$ and $\tilde{h}_{\alpha,n}$ are related to $G_n$.
We consider
$$C_n=\sqrt{n}\left(\|\bar{h}_{\alpha,n}-\bar{G}_n\|-\|\tilde{h}_\alpha-G\|\right)$$
and observe that $|A_n-C_n|\leq \sqrt{n}\|\tilde{h}_{\alpha,n}-\bar{h}_{\alpha,n}\|+\sqrt{n}\|G_n-\bar{G}_n\|=o_P(1)$ by (\ref{strongap}) (we are using 
that $\sqrt{n}\|\bar{U}_n-U_n \|\leq \|Q_n-\tilde{B}\|$, with  a similar
bound for the lower envelopes). Consequently, if suffices to prove convergence of $C_n$. From Theorem \ref{asymptotic1} we conclude that 
\begin{eqnarray*}
C_n \underset w \to \max\Big(\max_{t\in T_1}\tilde{B}(t),\max_{t\in T_2}(-\tilde{B}(t)),\max_{(s,t)\in T_3}{\textstyle \frac 1 2(\tilde{B}(t)-\tilde{B}(s))}\Big).
\end{eqnarray*}
The conclusion follows upon noting that (see Remark 4.2 in \cite{Hristo}) in the sets $T_i$, the function $G$ has local maxima:
if $t_0\in T_1$ then $G$ has a local maximum at $t_0$ and a local minimum if $t_0\in T_2$, also, if $(s_0,t_0)\in T_3$ then $G$ has a local maximum at $t_0$ and a local minimum 
at $s_0$. Therefore,
 $G'(t_0)=0$ and $G'(s_0)=0$ for every $t_0\in T_1,T_2$ or $(s_0,t_0)\in T_3$ and this entails
$w(t_0)=w(s_0)=\frac 1 {1-\alpha}$ for these points.
\hfill $\Box$

\end{document}